\documentclass{article}

\usepackage[a4paper,margin=1.7in]{geometry}

\usepackage{lipsum}
\usepackage{amsfonts}
\usepackage{amsfonts}
\usepackage{amsthm}
\usepackage{amsmath}
\usepackage{graphicx}
\usepackage{epstopdf}
\usepackage{algorithmicx}
\usepackage{algorithm}
\usepackage{algpseudocode}
\DeclareGraphicsExtensions{.eps,.pdf,.png,.jpg}
\usepackage{amsopn}

\title{Directional Total Generalized Variation Regularization}
\author{Rasmus Dalgas Kongskov\footnote{Department of Applied Mathematics and Computer Science, Technical University of Denmark, 2800 Kgs. Lyngby, Denmark (rara@dtu.dk, yido@dtu.dk, kiknu@dtu.dk .} \and Yiqiu Dong\footnotemark[1] \and Kim Knudsen\footnotemark[1]}

\usepackage{color}
\definecolor{pRed}{rgb}{0.7,0.01,0.15}
\definecolor{pBlue}{rgb}{0.01,0.01,0.8}
\definecolor{pGreen}{rgb}{0.01,0.8,0.01}
\definecolor{dGreen}{rgb}{0.01,0.4,0.01}

\usepackage[colorlinks=true,linkcolor=dGreen,citecolor=dGreen,urlcolor=dGreen]{hyperref}
\hypersetup{
  pdftitle={Directional Total Generalized Variation Regularization},
  pdfauthor={Rasmus Dalgas Kongskov, Yiqiu Dong and Kim Knudsen}
}

\graphicspath{{figs/}}
 
\newcommand{\qeda}{\hfill\ensuremath{\square}}

\newcommand{\dd}{\;\text{d}}

\newcommand{\figone}[3]{
\begin{figure}[h] 
\centerline{\includegraphics[width=#2\columnwidth]{#1}}
\caption{#3} \label{fig:#1}
\end{figure} }

\usepackage{cleveref}
\crefformat{equation}{(#2#1#3)}

\newcommand{\primf}{u} 				
\newcommand{\dom}{\Omega} 			
\newcommand{\TV}{\mathrm{TV}}			
\newcommand{\TGV}{\mathrm{TGV}}		
\newcommand{\DTGV}{\mathrm{DTGV}}		
\newcommand{\BGV}{\mathrm{BGV}}		
\newcommand{\BD}{\mathrm{BD}}			
\newcommand{\dualf}{v} 				
\newcommand{\dualfb}{\mathbf{\dualf}} 	
\newcommand{\scoord}{x} 				
\newcommand{\widE}{a}					
\newcommand{\angE}{\theta}			
\newcommand{\DTV}{\mathrm{DTV}}		
\newcommand{\noisef}{f}				
\newcommand{\regup}{\lambda}			
\newcommand{\scav}{t}					
\newcommand{\BV}{\mathrm{BV}}			
\newcommand{\intt}{n}					
\newcommand{\sysm}{A}					
\newcommand{\nv}{\mathbf{n}}			
\newcommand{\fo}{\sysm}				
\newcommand{\el}{E}					
\newcommand{\ball}{B_2(0)}			
\newcommand{\hball}{B_{2\times2}(0)}	
\newcommand{\rot}{R}					
\newcommand{\scm}{\Lambda}			
\newcommand{\efu}{\mathcal{J}}		
\newcommand{\con}{C}					
\newcommand{\eun}[1]{\left| #1 \right|_2}		
\newcommand{\xind}{i}					
\newcommand{\yind}{j}					
\newcommand{\fdm}{\nabla}				
\newcommand{\grad}{\nabla}			
\newcommand{\stdG}{\sigma}			
\newcommand{\kerG}{G_\stdG}			
\newcommand{\smp}{\mu}				
\newcommand{\aone}{\primf}			
\newcommand{\atwo}{\dualf}			
\newcommand{\cosm}{c}					
\newcommand{\sinm}{s}					
\newcommand{\odiff}{e}				
\newcommand{\niter}{k}				
\newcommand{\imq}{q}					
\newcommand{\dualCP}{p}				
\newcommand{\stepA}{\eta}				
\newcommand{\stepB}{\tau}				
\newcommand{\tolr}{\mathrm{tol}}		
\newcommand{\lipc}{L}					
\newcommand{\dualt}{\mathbf{w}}		
\newcommand{\dualtm}{W}				
\newcommand{\dualte}{w}				
\newcommand{\kTGV}{h}					
\newcommand{\lTGV}{l}					
\newcommand{\ddiv}{\widetilde{\mathrm{div}}} 
\newcommand{\kset}{K}					
\newcommand{\projP}{P}				
\newcommand{\affP}{\mathcal{A}^1}		
\newcommand{\projS}{\mathcal{S}}		
\newcommand{\dide}{D}					
\newcommand{\forw}{+}					
\newcommand{\bacw}{-}					
\newcommand{\funcB}{\mathcal{F}}		
\newcommand{\noise}{\eta}				
\newcommand{\pdim}{M}					
\newcommand{\inddiv}{b}				
\newcommand{\inddivB}{\gamma}			
\newcommand{\Mset}{N}					
\newcommand{\rchi}{\raisebox{2pt}{$\chi$}}       

\theoremstyle{plain}
\newtheorem{thm}{Theorem}[section] 

\theoremstyle{definition}
\newtheorem{defn}[thm]{Definition} 	
\newtheorem{exmp}[thm]{Example} 		
\newtheorem{prop}{Proposition}		
\newtheorem{cor}{Corollary}			
\newtheorem{remark}[thm]{Remark} 

\begin{document}

\maketitle
 
\begin{abstract}
In inverse problems, prior information and a priori-based regularization techniques play important roles. In this paper, we focus on image restoration problems, especially on restoring images whose texture mainly follow one direction. In order to incorporate the directional information, we propose a new directional total generalized variation (DTGV) functional, which is based on total generalized variation (TGV) by Bredies \textit{et al}. [SIAM J. Imaging Sci., 3 (2010)]. After studying the mathematical properties of DTGV, we utilize it as regularizer and propose the L$^2$-$\DTGV$ variational model for solving image restoration problems. Due to the requirement of the directional information in DTGV, we give a direction estimation algorithm, and then apply a primal-dual algorithm to solve the minimization problem. Experimental results show the effectiveness of the proposed method for restoring the directional images. In comparison with isotropic regularizers like total variation and TGV, the improvement of texture preservation and noise removal is significant. 
\end{abstract}

\section{Introduction}

In the field of inverse problems, regularization techniques have been introduced to overcome the ill-posedness in order to obtain reasonable and stable solutions. For many image processing problems incorporating prior information through regularization techniques has attracted much attention. In this paper, we will study directional regularization for image restoration problems.

The image is given in the domain $\dom$, a connected bounded open subset of $\mathbb{R}^2$ with Lipschitz boundary, and given by a real-valued function  $\hat{\primf}: \dom\rightarrow\mathbb{R}$. The image is degraded through an operator $\fo\in\mathcal{L}(L^2(\dom))$ and by additive white Gaussian noise $\noise$, and thus the degraded image $\noisef$ is given by 
\begin{align} \label{eqn:forward}
\noisef = \fo\hat{\primf} +  \noise.
\end{align}
We consider $A$ as the identity operator (denoising problem) and $A$ having the form of a blurring operator (deblurring problem). 
The analysis is based on the variational model for image restoration
\begin{align}\label{eqn:general}
\min_\primf\ \frac12 \|\fo u-\noisef\|^2_{L^2(\dom)} + \regup \mathcal{R}(\primf),
\end{align}
where $\mathcal{R}$ is the regularization term, which incorporates prior information on $\hat{\primf}$, and $\regup>0$ is the regularization parameter, which controls the trade-off between the fit with the data $\noisef$ and the regularization. 

Due to its capability of preserving sharp edges, total variation (TV) regularization proposed in \cite{Rudin1992} has been used for many image processing problems, e.g.\ in image denoising \cite{Vogel1996, dong2009efficient, dong2013convex, sciacchitano2015variational}, in blind deconvolution \cite{Chan1998}, in tomographic reconstruction \cite{Delaney1998, Jonsson}, etc. Although TV regularization is very effective for restoring piece-wise constant images, it has some shortcomings, and the most notable one is the appearance of staircasing artifacts in slanted regions \cite{Nikolova2000,Ring2000}. To overcome staircasing artifacts, higher-order derivatives have been used, see \cite{Scherzer1998,Chan2000,Setzer2008b,Setzer2011}. In \cite{Bredies2010}, total generalized variation (TGV) of order $h$, was proposed, which incorporates the first up to the $h$-th derivatives. When $h$ equals 1, it yields the TV regularization. 

\begin{figure}[t]
  \centering
	\begin{minipage}[b]{0.27\textwidth}
	\includegraphics[width=\textwidth]{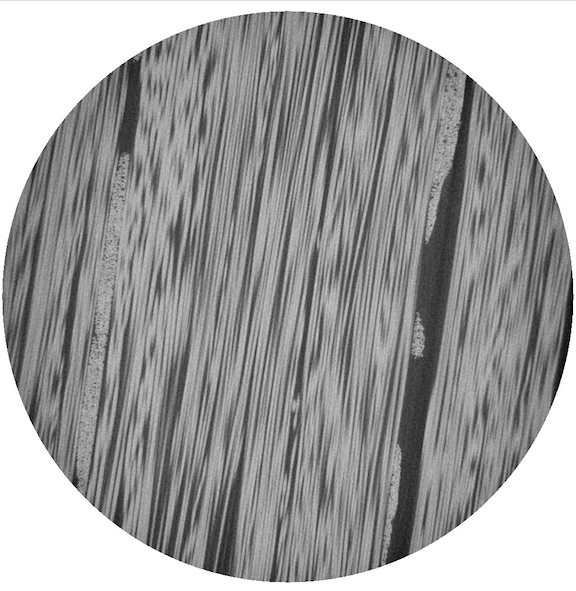}
	\end{minipage}
\hspace{1.5cm}
	\begin{minipage}[b]{0.57\textwidth}
	\includegraphics[width=\textwidth]{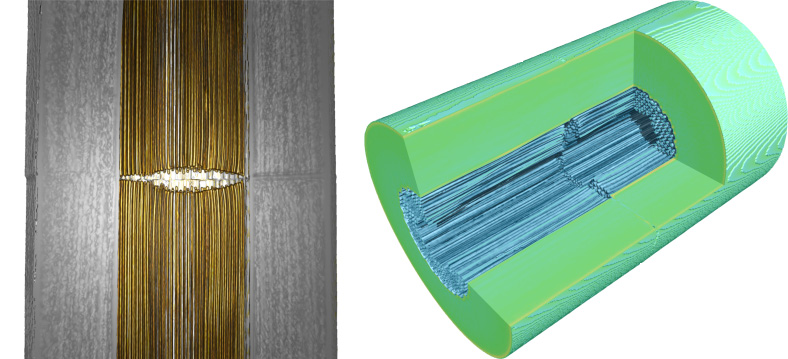}
	\end{minipage}
\caption{Left: A CT scan of uni-directional glass fibre (obtained from \cite{Jespersen2016b}, see more in \cite{Jespersen2016}). Right: A CT scan and 3D-model of an optical fibre with a cavity (obtained from \cite{Sandoghchi2014}).} \label{fibre}
\end{figure}

In many applications related to fibers, the textures in images have very clear directionality. Examples include glass fibres in wind-turbine blades, optical fibres for communication, and ceramic fibres in fuel cells; see \cref{fibre}. Another application with clear directional textures is seismic imaging. We call images with textures oriented mainly along one certain direction {\it directional images}. 
Achieving high-quality images is crucial for the analysis of these fibre materials, therefore imposing the directional information of the texture is highly desirable. 

Directional regularization has been introduced for standard TV in \cite{Esedoglu2004,Berkels2006,Turgay2009, Fei2013,Fernandez-Granda2013} and in terms of shearlet-based TV in \cite{Easley2009}. In \cite{Bayram2012} a type of directional TV is introduced for image denoising based on images with one main direction. This method is further developed to be spatially adaptive in \cite{Zhang2013} via pixel-specific angle-estimates. Moreover, the directional information through the structure tensor, defined in \cite{Weickert1998a,Weickert1999a}, has been used to extend TV regularization. This new regularization method is called structure tensor total variation (STV), which has been applied in different imaging problems, see \cite{Lefkimmiatis2015, Estellers2015}. The structure tensor has also been combined with second-order-derivative regularization in \cite{Hafner2015}. In \cite{Ranftl2012,Ferstl2013}, anisotropic diffusion tensor has been applied only on the first-order-derivative term in the second-order TGV. All these directional regularization techniques are introduced for discretized problems, but the underlying continuous problems are not studied. Furthermore, it is in these previous work not clear if directional information can be incorportaed in higher-order derivatives, e.g. through TGV regularizer.

The first goal of this paper is to formulate in a continuous setting regularization terms that incorporate directionality.  We first define the directional TV (DTV) functional and generalize it to higher orders, the so-called directional TGV (DTGV). We construct this generalization to higher orders such that the directional information is also incorporated in higher-order derivatives, which is different from the anisotropic TGV proposed in \cite{Ranftl2012,Ferstl2013}. Under the continuous setting we study the mathematical properties of the DTGV functional. Further, we utilize the DTGV functional as a regularizer in \cref{eqn:general} and derive existence and uniqueness results for the minimization problem in \cref{eqn:general}. The second goal of the paper is to give a numerical implementation based on the primal-dual algorithm proposed in \cite{Chambolle2011a} to solve the minimization problem in \cref{eqn:general}, and through numerical experiments evaluate its performance. Since DTGV requires the input of the main direction we also propose a direction estimation algorithm.



The rest of the paper is organized as follows. In \cref{sec:DTV} we define the directional total variation (DTV) functional. Through two equivalent definitions of DTV, we obtain a hint of how to incorporate directional information into TGV. In \cref{sec:DTGV}, we propose the second order directional total generalized variation ($\DTGV^2_\regup$) functional, and extend it to higher orders. We study the mathematical properties of DTGV in \cref{sec:properties}, and in \cref{sec:rergularization} we apply it as regularization in \cref{eqn:general} to propose a new L$^2$-$\DTGV^2_\regup$ model. The existence and uniqueness results for the L$^2$-DTGV$^2_\regup$ model is also provided. In \cref{sec:algorithms} we introduce a direction estimation algorithm in order to obtain the required main direction from the degraded images, and then apply a primal-dual algorithm for solving the minimization problem in our restoration model based on the work proposed in \cite{Chambolle2011a}. The numerical results shown in \cref{sec:numexp} demonstrate the effectiveness of the direction estimation algorithm, the influence of the parameters in DTGV, and the performance of our restoration method. Finally, conclusions are drawn in \cref{sec:conclusion}.

\section{Directional Total Variation} \label{sec:DTV}

The definition of total variation (TV) for $\primf \in \BV(\dom)$, the space of functions of bounded variation over the domain $\dom$, can be written as \cite{Aubert2006,Rudin1992}
\begin{align} \label{eqn:TVform}
\TV(\primf)
&= \sup \left\{ \left. \int_\dom \primf \,\, \mathrm{div}\, \dualfb \,\dd \scoord \right| \dualfb\in C^1_c(\dom,\mathbb{R}^2), \dualfb(\scoord) \in \ball \,\, \forall \scoord \in \dom  \right\},
\end{align}
where $\dualfb$ denotes the dual-variable and $\ball$ denotes the closed Euclidean unit ball centered at the origin. In this section, we will introduce directional information into TV and define directional total variational (DTV). The idea of DTV was first proposed in \cite{Bayram2012} in the discrete case. Following a similar idea we will give the DTV definition in the continuous case. Through  examples we demonstrate the differences between TV and DTV. More mathematical properties will be derived based on the extension to total generalized variation (TGV) in \cref{sec:DTGV}. 

TV is rotational invariant. In order to allow rotational variation, we restrict the dual variable $\dualfb$ in an ellipse instead of the unit ball. 
Define the closed elliptical set, $\el^{\widE,\angE}(0)$, centered at the origin with the major semi-axis 1 oriented in direction $(\cos\angE,\sin\angE)$ and the minor semi-axis $\widE \in (0,1]$ by
\begin{align}\label{def:Ellip}
\el^{\widE,\angE}(0) = \left\{ \left(\begin{array}{c}x_1 \\ x_2 \end{array}\right)\in \mathbb{R}^2 \left| \left( \frac{x_1\cos\angE + x_2\sin\angE}{1} \right)^2 + \left( \frac{-x_1\sin\angE + x_2\cos\angE}{a}  \right)^2 \leq 1 \right. \right\}.
\end{align}
In \cref{fig:ellipse2} the elliptical set $\el^{\widE,\angE}(0)$ is depicted. We are now ready to define DTV:
\figone{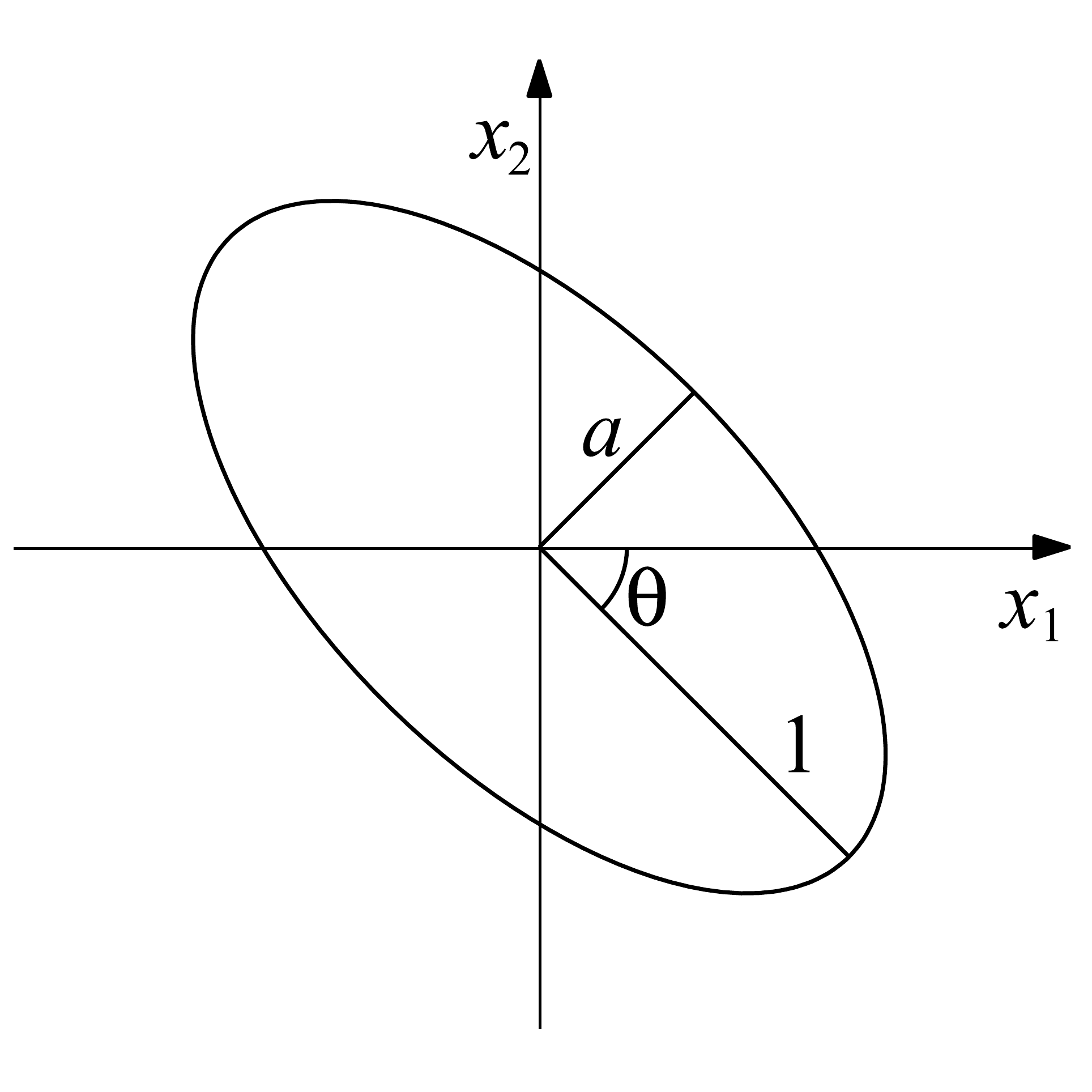}{0.27}{Sketch of the elliptical set $\el^{\widE,\angE}(0)$. Here it is shown for $\widE=0.5$ and $\angE=-\frac{\pi}{4}$.}
\begin{defn} \label{def:DTV1} The directional total variation $(\mathrm{DTV})$ with respect to $(\widE,\angE)$ of a function $\primf\in\BV(\dom)$ is defined as
\begin{align*} 
\DTV(\primf)
&= \sup \left\{ \left. \int_\dom \primf \,\, \mathrm{div}\, \tilde{\dualfb} \,\dd \scoord\right| \tilde{\dualfb} \in C^1_c(\dom,\mathbb{R}^2), \tilde{\dualfb}(\scoord) \in \el^{\widE,\angE}(0) \,\, \forall \scoord \in \dom  \right\}.
\end{align*}
\end{defn}


Introduce the rotation matrix $\rot_\angE$ and the translation matrix $\scm_\widE$ by
\begin{align*} 
\rot_\angE = \begin{pmatrix} \cos\angE & -\sin\angE \\ \sin\angE & \cos\angE \end{pmatrix} \mbox{ and} \quad \scm_\widE = \begin{pmatrix} 1 & 0 \\ 0 & \widE \end{pmatrix}.
\end{align*}
Then 
\begin{align} \label{eqn:dual1to2}
\tilde{\dualfb}(x) = \rot_\angE \scm_\widE \dualfb(x) \in \el^{\widE,\angE}(0) \qquad \Leftrightarrow \qquad \dualfb(x)= \scm_{\frac{1}{\widE}} \rot_{-\angE} \tilde{\dualfb}(x)\in \ball. 
\end{align}
Further, we define the directional divergence for $\dualfb(x)\in\ball$
\begin{equation}\label{ddiv1d}
\ddiv \dualfb(x)=\mathrm{div}\,  \rot_\angE\scm_\widE \dualfb(x)=\mathrm{div}\tilde{\dualfb}(x).
\end{equation}
Using the relations in \cref{eqn:dual1to2} and the definition in \cref{ddiv1d} we give another equivalent definition for DTV.
\begin{defn} \label{def:DTV2}For a function $\primf\in\BV(\dom)$, its DTV is defined as
\begin{align*} 
\DTV(\primf)
&= \sup \left\{ \left. \int_\dom \primf \,\, \ddiv \dualfb  \,\dd \scoord \right| \dualfb\in C^1_c(\dom,\mathbb{R}^2), \dualfb(\scoord) \in \ball \,\, \forall \scoord \in \dom  \right\}.
\end{align*}
\end{defn}

This definition is very similar to the one for TV in \cref{eqn:TVform}; the only difference is the change on the divergence operator in the integral. Here we provide three examples on comparisons of DTV with TV:
\begin{exmp}
  For $u\in C^\infty_c(\Omega)$ we have using integration by parts
  \begin{align*}
    \TV(u) &= \|\nabla u \|_{L^1(\Omega)},\\
    \DTV(u ) &= \left \|
           \begin{pmatrix}
             D_\theta u\\
             a D_{\theta^\perp} u, 
           \end{pmatrix}\right \|_{L^1(\Omega)}
  \end{align*}
  where $D_\theta u$ denotes the directional derivative of $u$ in the direction $(\cos\theta,\sin \theta)$ and $\theta^\perp = \theta + \pi/2.$ 
\end{exmp}
The above example shows that $\DTV(u)$ is an anisotropic total variation functional. The next example demonstrates the difference between DTV and TV.  
\begin{exmp} 
Define $\primf_1=\rchi_{\ball}$, the characteristic function of the unit disk. 
The total variation of such a characteristic function is given by the length of the perimeter, i.e. $\TV(\primf_1) = 2\pi$. The calculation of $\DTV(\primf_1)$ is according to \cref{def:DTV2}: by using the divergence theorem with $\nv$ denoting the outward unit normal vector we have for any $\dualfb\in C_c^1(\Omega,\mathbb R^2)$
\begin{align*}
\int_\dom \primf_1 \,\, \ddiv \dualfb  \,\dd \scoord &=\int_{{\ball}} \mathrm{div} \rot_\angE\scm_\widE \dualfb \dd \scoord\\ 
&= \int_{\partial {\ball}} \left( \rot_\angE \scm_\widE \dualfb \right) \cdot \nv \dd s \\
&=  \int_{\partial{\ball}} \dualfb \cdot  \left( \scm_\widE\rot_{-\angE} \nv \right) \dd s. 
\end{align*}
The integrand is maximized among unit vector fields for $\dualfb =  \left(  \scm_\widE\rot_{-\angE}  \nv \right)/ |  \scm_\widE\rot_{-\angE}  \nv |$ thus yielding
\begin{align}\label{dtvu1}
   \DTV(u_1) &= \int_{\partial {\ball}} | \scm_\widE\rot_{-\angE}  \nv | ds \\
   &= \int_0^{2\pi} (\cos^2 (\tau-\angE) + a^2 \sin^2(\tau-\angE))^{1/2} d\tau \notag\\
   &=\int_0^{2\pi} (\cos^2 \omega + a^2 \sin^2 \omega)^{1/2} d\omega, \notag
\end{align}
which is by the way the length of the perimeter of the elliptical set $\el^{a,\theta}(0).$ 
\end{exmp}

\begin{exmp}
Take now instead the characteristic function $\primf_2=\rchi_{p\el^{b,0}(0)}$ of the elliptical set $p\el^{b,0}(0)$ with $0<b<1.$ By choosing $p$ such that
\begin{align*}
  p {\int_0^{2\pi} (b^{2}\cos^2 \tau + \sin^2\tau)^{1/2} d\tau}  = {2\pi}
\end{align*}
the length of the ellipse perimeter is $TV(u_2) = 2\pi$ as before. 

Let us compute $\DTV(u_2)$ for the  two different orientations given by $\theta = 0$ and $\theta=\pi/2.$ First fix $\theta = 0$ and $0<a< 1.$ The outward unit normal to $p\el^{b,0}(0)$ at point $(p\cos\tau,pb\sin\tau)$ is
\begin{align}\label{normal}
  n(\tau)  = \frac 1 {(b^2\cos^2\tau + \sin^2\tau)^{1/2}}
  \begin{pmatrix}
    b\cos \tau\\ \sin\tau
  \end{pmatrix}
\end{align}
and then with $\angE=0$ calculate 
\begin{align}\label{angle0}
  \int_{\partial{p\el^{b,0}(0)}} \left( \rot_\angE \scm_\widE \dualfb \right) \cdot \nv \dd s &= \int_0^{2\pi}(\dualfb \cdot \scm_{\widE}\rot_{-\angE}\nv)|x'(\tau)|\ d\tau\\
&= p\int_0^{2\pi} \left( \dualfb \cdot \begin{pmatrix}
    b  \cos \tau\\ a \sin\tau
  \end{pmatrix}\right)\ d\tau.\notag
\end{align}
The integrand is maximized with the unit vector field
\begin{align*}
  \dualfb(\tau) = \frac 1 {(b^2\cos^2\tau + a^2\sin^2\tau)^{1/2}} \begin{pmatrix}
     b \cos \tau\\ a \sin\tau
  \end{pmatrix}
\end{align*}
and thus we obtain
\begin{align}\label{dtvu2}
  d_0= \DTV(u_2) &= p \int_0^{2\pi} \left( b^2\cos^2\tau + a^2\sin^2\tau\right)^{1/2} \ d\tau
 \end{align}
By using Maple to estimate the integration numerically, we compare \eqref{dtvu1} and \eqref{dtvu2}, and find that while $\TV(u_1) = \TV(u_2),$ $\DTV(u_2) < \DTV(u_1).$

Next, we compute $\DTV(u_{2})$ with $\angE=\frac\pi2$ and $0<\widE\leq1$. In this case, according to \cref{angle0} by using \cref{normal} we get
\[
  \int_{\partial{p\el^{b,0}(0)}} \left( \rot_\angE \scm_\widE \dualfb \right) \cdot \nv \dd s = p\int_0^{2\pi} \left( \dualfb \cdot \begin{pmatrix}
   -\sin \tau\\ \widE b \cos\tau
  \end{pmatrix}\right)\ d\tau.
\]
Its maximum is reached at 
\begin{align*}
  \dualfb(\tau) = \frac 1 {(\sin^2\tau + a^2b^{2}\cos^2\tau)^{1/2}} \begin{pmatrix}
     -\sin \tau\\ ab \cos\tau
  \end{pmatrix}
\end{align*}
and we obtain
\begin{align}\label{dtvu3}
  d_{\pi/2} = \DTV(u_2) &= p \int_0^{2\pi} \left( \sin^2\tau + a^2b^{2}\cos^2\tau\right)^{1/2} \ d\tau
 \end{align}
By using Maple to estimate the integration numerically, we compare \eqref{dtvu2} and \eqref{dtvu3}, and find $d_0 \leq d_{\pi/2}.$ 

\begin{figure}[t] 
\centerline{\includegraphics[width=0.4\textwidth]{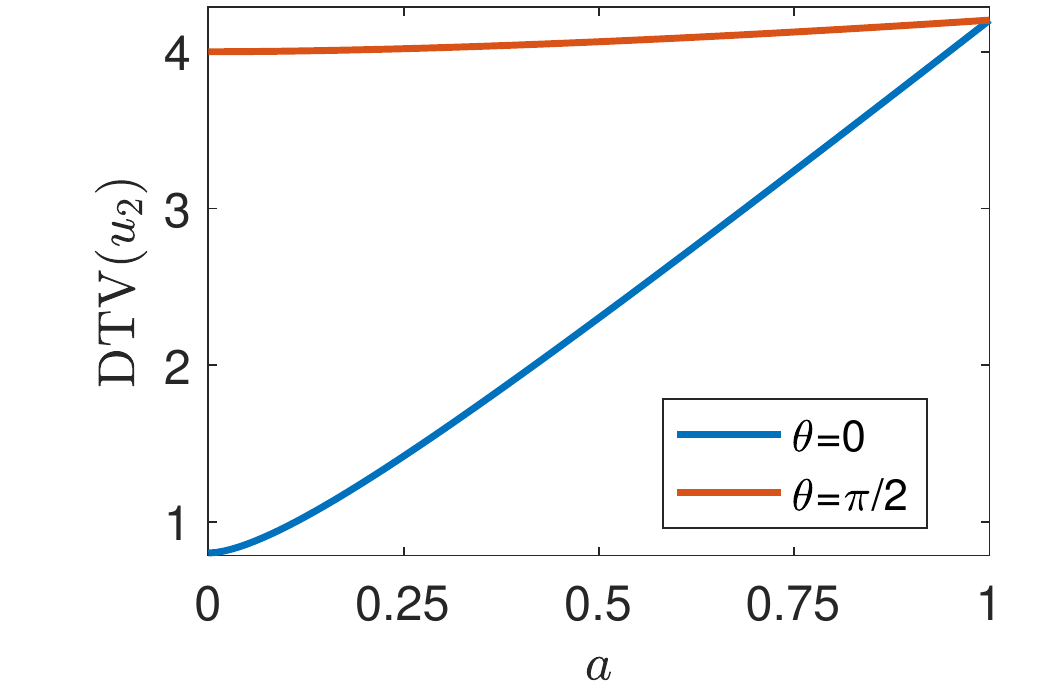} \hspace{10mm}
\includegraphics[width=0.4\textwidth]{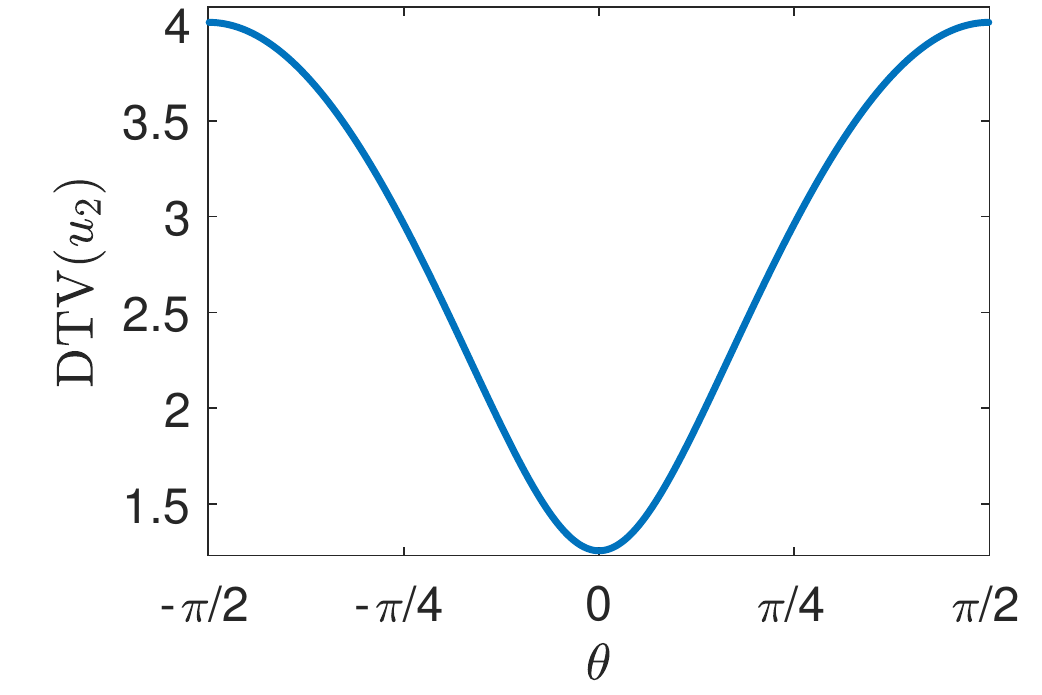}}
\caption{For $\primf_2$ with $b = 0.2.$ Left: $\DTV(u_{2})$ as a function of $a$ and three choices of $\angE$. Right: $\DTV(u_{2})$ as a function of $\theta$ and $\widE=0.5.$} \label{fig:example1}
\end{figure}

In \cref{fig:example1} we illustrate the dependency of DTV on $\angE$ and $\widE$.  The left plot shows that for fixed $\angE,$ $\DTV(u_2)$ is monotonically incresing with $\widE.$ Moreover, independent of the value of $\widE,$ a smaller and more correct choice of $\angE$, i.e. closer to the main direction of the object, gives a lower $\DTV(u_2).$ However, for $\widE=1$ the values of $\DTV(u_2)$ agree, since in that case DTV is equivalent to the rotationally invariant TV. The right plot suggests that $\DTV(u_2)$ depends as a scaled and translated sine function with a minimum obtained when the rotation angle $\angE$ in DTV coincides with the main direction of the object (i.e.\ $\angE=0$). 


\end{exmp}

As a side remark we note that the extension of TV to DTV is based on knowledge (or estimates of) one global direction. Extensions of the presented work to the important case of spatially dependent directions is left for future work.


In imaging problems a common artifact caused by TV regularization is staircasing, which is a classical example of a mismatch between prior knowledge and the reality, i.e., smooth regions are approximated by piece-wise constant regions. One way to overcome the staircasing effect is to use higher order derivatives in the regularization; this is the topic of the next section.

\section{Directional Total Generalized Variation} \label{sec:DTGV}
Total Generalized Variation (TGV) is a generalization of TV to a functional, which takes derivatives of order $\kTGV>0$ into account, and the first order TGV, i.e., $\kTGV=1$, is identical to TV. 
It turns out that for natural images, TGV regularization for denoising is often superior to TV regularization, and the staircasing effect is well avoided \cite{Bredies2010}. 

For a 2-by-2 symmetric matrix-valued function  
$$V(x) = \begin{pmatrix}v_{11}(x) & v_{12}(x) \\ v_{12}(x) & v_{22}(x)  \end{pmatrix} \in \mathrm{Sym}^2(\mathbb{R}^2) $$ define
\begin{align*}
 (\mathrm{div}\, V)^\top = \left(\begin{matrix}\frac{\partial v_{11}}{\partial \scoord_1} + \frac{\partial v_{12}}{\partial \scoord_2}\\ \frac{\partial v_{21}}{\partial \scoord_1} + \frac{\partial v_{22}}{\partial \scoord_2}\end{matrix}\right), \qquad 
&\mathrm{div}^2\, V = \frac{\partial^2 v_{11}}{\partial \scoord_1^2} + \frac{\partial^2 v_{22}}{\partial \scoord_2^2} + 2\frac{\partial^2 v_{12}}{\partial \scoord_1 \partial \scoord_2}.
\end{align*}
Define further the second order unit sphere $\hball$ consisting of
matrices $V = (v_{i,j})_{i,j = 1}^2$ with rows $(v_{i1}, v_{i2})^\top\in\ball$ and
columns $(v_{1i}, v_{2i})^\top\in\ball$ for $i=1, 2$. We write $W \in \regup_0\hball$ if $ W/\lambda_0 \in \hball$ with $\lambda_0 >0$, leaving $\|W\|_2 \leq \lambda_0$. Then we can for a
function $\primf \in L^1(\Omega)$ define the second order TGV for $\lambda=(\lambda_0,\lambda_1)$ by 
\begin{align}\label{def:TGV2}
\TGV^2_\regup(\primf) = \sup& \left\{ \left. \int_\dom \primf \,\, \mathrm{div}^2\, \dualtm\,\dd \scoord \right|\right. \\
&\left.\dualtm\in C^2_c(\dom,\mathrm{Sym}^2(\mathbb{R}^2)), \dualtm(x)\in \regup_0\hball, (\mathrm{div}\, \dualtm(x))^\top \in \regup_1\ball\, \forall x\in\dom\right\}. \notag
\end{align}
 (Note that in \cite{Bredies2010} higher order TGV is defined equivalently using the Frobenius norm.)   

In order to include directional information in $\TGV^{2}_{\regup}$ we replace higher order balls by higher order elliptical sets. Denote by $\el_{2\times2}^{\widE,\angE}(0)$ the space of matrices  $V = (v_{i,j})_{i,j = 1}^2$ with rows $(v_{i1}, v_{i2})^\top\in \el^{\widE,\angE}(0) $ and
columns $(v_{1i}, v_{2i})^\top\in\el^{\widE,\angE}(0)$ ($i=1, 2$).  
Then we can define the second order directional TGV (DTGV):
\begin{defn} The second order directional total generalized variation, $\DTGV^2_\regup$, of a function $\primf\in L^1(\dom)$ is defined as
\begin{align} \label{def:DTGV2a}
\DTGV^2_\regup(\primf) = &\sup \left\{ \left. \int_\dom \primf \,\, \mathrm{div}^2\, \widetilde{\dualtm}\,\dd \scoord \right| \right.\\
&\left.\widetilde{\dualtm} \in C^2_c(\dom,\mathrm{Sym}^2(\mathbb{R}^2)), \widetilde{\dualtm}(x)\in \regup_0\el_{2\times2}^{\widE,\angE}(0), (\mathrm{div}\, \widetilde{\dualtm}(x))^\top \in \regup_1\el_{2}^{\widE,\angE}(0)\,\, \forall x\in\dom \right\}.\notag
\end{align}
\end{defn}

\begin{remark}
In \cite{Holler2014}, TGV has been extended to the infimal convolution of a number of TGV type functionals with arbitrary norms. Our definition of $\DTGV^2_\regup$ in \eqref{def:DTGV2a} can be considered as a special case under this general definition  with norm defined by particular ellipses. 
\end{remark}

Next, we will provide a characterization of $\DTGV^{2}_{\regup}$, which has the same feasible set as in \cref{def:TGV2}. 
To do so we need the second order generalization of \eqref{eqn:dual1to2} 
\begin{align} \label{eqn:tildeWaW}
\widetilde{V}(x) = \rot_\angE\scm_\widE V(x) \scm_\widE\rot_{-\angE} \in \el_{2\times2}^{\widE,\angE}(0) \qquad \Leftrightarrow \qquad V(x) \in \hball,
\end{align}
which follows easily from the definitions.
In addition, we define for $W\in C_c^2(\Omega,\mathrm{Sym}^{2}(\mathbb{R}^{2}))$ the directional divergence
\begin{align} \label{eqn:ddiv}
\ddiv W(x) = \mathrm{div} \rot_\angE\scm_\widE W(x)
\end{align}
and the second order directional divergence
\begin{equation} \label{eqn:ddiv2}
\ddiv^2 W(x) = \mathrm{div}^2 \, \widetilde{\dualtm}(x) \qquad \mbox{with} \qquad \widetilde{\dualtm}(x) = \rot_\angE\scm_\widE\dualtm(x)\scm_\widE\rot_{-\angE}.
%
\end{equation}
The characterization of  $\DTGV^{2}_{\regup}$ is now as follows:
\begin{thm} \label{thm:DTGVdef}
With the directional divergence and the second order directional divergence defined in \cref{eqn:ddiv} and \cref{eqn:ddiv2}, for  
$\primf\in L^1(\dom)$
\begin{align}\label{def:DTGV2b}
\DTGV^2_\regup(\primf) = \sup &\left\{ \left. \int_\dom \primf \,\, \ddiv^2\, \dualtm\,\dd \scoord \right|\right.\\
 &\left.\dualtm\in C^2_c(\dom,\mathrm{Sym}^2(\mathbb{R}^2)), \|\dualtm(x)\|_2 \leq \regup_0,\|(\ddiv\, \dualtm(x))^{\top}\|_2 \leq \regup_1\,\, \forall x\in\dom\right\}.\notag
\end{align}
\end{thm}
\textit{Proof:} 
Let $\widetilde{\dualtm}  \in C^2_c(\dom,\mathrm{Sym}^2(\mathbb{R}^2))$ with $\widetilde{\dualtm}(x)\in \regup_0\el_{2\times2}^{\widE,\angE}(0)$ and $(\mathrm{div}\, \widetilde{\dualtm}(x))^\top \in \regup_1\el_{2}^{\widE,\angE}(0)$ for any $x\in\dom$. Due to \cref{eqn:tildeWaW} we have for any $x\in\dom$
\[
\dualtm(x)=\scm_{\frac{1}{\widE}}\rot_{-\angE}\widetilde{\dualtm}(x)\rot_{\angE}\scm_{\frac{1}{\widE}} \in \regup_{0}\hball,
\] 
and hence $\|W(x)\|_2\leq \lambda_0.$ In addition, since $(\mathrm{div}\, \widetilde{\dualtm}(x))^\top \in \regup_1\el_{2}^{\widE,\angE}(0)$, we obtain  according to  \cref{eqn:dual1to2} that
\[
\scm_{\frac{1}{\widE}}\rot_{-\angE}(\mathrm{div}\, \widetilde{\dualtm}(x))^\top=\scm_{\frac{1}{\widE}}\rot_{-\angE}(\mathrm{div}\, \rot_\angE\scm_\widE \dualtm(x) \scm_\widE\rot_{-\angE})^{\top}=(\ddiv\,\dualtm(x))^{\top}\in\regup_{1}\ball,
\]
i.e.\ $\|(\ddiv\, \dualtm(x))^{\top}\|_2 \leq \regup_1$. Hence, we have proven that the feasible set in \cref{def:DTGV2b} is equivalent to the one in \cref{def:DTGV2a}, and the result follows from \eqref{eqn:ddiv2}.

\qeda 

To close the section we  extend the definition of $\DTGV^2_{\regup}$  to arbitrary order $h\in\mathbb{N}.$  Recall from \cite{Bredies2010} the higher order TGV defined for any $\kTGV\in\mathbb{N},\; \lambda = (\lambda_0,\lambda_1,\ldots,\lambda_{h-1})$ by
\begin{align}\label{eqn:TGV}
\TGV^\kTGV_\regup(\primf) = \sup &\left\{ \left. \int_\dom \primf \,\, \mathrm{div}^\kTGV \dualt\,\dd \scoord \right|\right.\\
 &\left.\dualt\in C^\kTGV_c(\dom,\mathrm{Sym}^\kTGV(\mathbb{R}^2)), \|{\mathrm{div}^\lTGV\, \dualt(x)}\|_2 \leq \regup_\lTGV,\,\, \forall x\in\dom \,\,\mbox{and}\,\, \lTGV=0,...,\kTGV-1 \right\},\notag
\end{align} 
where $\mathrm{Sym}^\kTGV(\mathbb{R}^2)$ denotes the space of symmetric $\kTGV$-tensors in $\mathbb{R}^2$. For any $x\in\dom$, $\dualt(x)$ is a symmetric $\kTGV$-tensor. The operator $\mathrm{div}^\lTGV$ on $\dualt(x)$ is defined as
\[
(\mathrm{div}^\lTGV \dualt(x))_\inddiv\, = \sum_{\inddivB\in \Mset_\lTGV} \frac{\lTGV!}{\inddivB!} \frac{\partial^\lTGV \dualt(x)_{\inddiv+\inddivB}}{\partial \scoord^\inddivB} \quad \text{for each component} \quad \inddiv\in\Mset_{\kTGV-\lTGV}
\]
where $\Mset_\kTGV = \left\{\inddiv\in\mathbb{N}^2 \left| |\inddiv| = \sum_{i=1}^2 \inddiv_i = \kTGV  \right.\right\}$. In addition, $\|\cdot\|_2$ on a symmetric $\kTGV$-tensor is defined as
\begin{align*}
\|\dualt(x)\|_2 = \left(\sum_{\inddiv \in \Mset_\kTGV} \frac{\kTGV!}{\inddiv!} (\dualt(\scoord))_\inddiv^2 \right)^{1/2}.
\end{align*}
For further details of the $\kTGV$-tensors we refer the readers to \cite{Bredies2010}.

To take the directional information into account, we first define two tensor fields:
\begin{align}
B_{\footnotesize\underbrace{2\times \ldots\times 2}_{m \text{ times}}}(0)=&\left\{\xi\in\mathrm{Sym}^{m}(\mathbb{R}^{2}): \xi(e_{p_{1}},\cdots, e_{p_{i-1}},\cdot, e_{p_{i+1}},\cdots, e_{p_{m}})\in B_{2}(0)\right.\label{high_orderB}\\
&\left.  \mbox{for any } i\in\{1,\cdots,m\} \mbox{ and } p\in\{1, 2\}^{m}\right\},\notag\\
\el^{\widE,\angE}_{\footnotesize\underbrace{2\times \ldots\times 2}_{m \text{ times}} }(0)=&\left\{\xi\in\mathrm{Sym}^{m}(\mathbb{R}^{2}): \xi(e_{p_{1}},\cdots, e_{p_{i-1}},\cdot, e_{p_{i+1}},\cdots, e_{p_{m}})\in E^{\widE, \angE}_{2}(0)\right.\label{high_orderE}\\
&\left.  \mbox{for any } i\in\{1,\cdots,m\} \mbox{ and } p\in\{1, 2\}^{m}\right\},\notag
\end{align}
where $e_{1}$ and $e_{2}$ denote the two standard basis vectors in $\mathbb{R}^{2}$. Similar as in \cref{ddiv1d} and \cref{eqn:ddiv2} we define $\lTGV$-directional divergence, $\ddiv^\lTGV$ for any order $\lTGV\leq\kTGV$ as
\[
\ddiv^{l}\dualt(x)=\mathrm{div}^l \widetilde{\dualt}(x),
\]
where $\dualt(x)\in B_{\footnotesize\underbrace{2\times \ldots\times 2}_{h \text{ times}}}(0)$, and we can obtain $\widetilde{\dualt}(x)\in \el^{\widE,\angE}_{\footnotesize\underbrace{2\times \ldots\times 2}_{l \text{ times}} }(0)\otimes B_{\footnotesize\underbrace{2\times \ldots\times 2}_{h-l \text{ times}}}(0)$ referring to the relations in \cref{eqn:dual1to2} with $\otimes$ as the tensor product. Then as in \cref{def:DTGV2a}, we give the definition of $\DTGV^\kTGV_\regup$ as follows:
\begin{defn} The $\kTGV$'th-order directional total generalized variation, $\DTGV^\kTGV_\regup$, of a function $\primf\in L^1(\dom)$ is defined as
\begin{align} \label{def:DTGV2}
\DTGV^\kTGV_\regup(\primf) &= \sup \bigg\{ \left. \int_\dom \primf \,\, \mathrm{div}^\kTGV \dualt\,\dd \scoord \right|\\
 &\dualt\in C^\kTGV_c(\dom,\mathrm{Sym}^\kTGV(\mathbb{R}^2)), \mathrm{div}^\lTGV\,\dualt(x) \in \regup_l\el^{\widE,\angE}_{\footnotesize\underbrace{2\times \ldots\times 2}_{\kTGV-\lTGV\text{ times}} }(0),\,\, \forall x\in\dom \,\,\mbox{and}\,\, \lTGV=0,...,\kTGV-1 \bigg\}.\notag\end{align}
\end{defn}

Following the similar idea as the proof of \cref{thm:DTGVdef}, we obtain a characterization of $\DTGV^\kTGV_\regup$.
\begin{thm} For $\primf\in L^1(\dom)$ 
\begin{align} \label{def:DTGV}
\DTGV^\kTGV_\regup(\primf) = \sup &\left\{ \left. \int_\dom \primf \,\, \ddiv^\kTGV \dualt\,\dd \scoord \right|\right.\\
 &\left.\dualt\in C^\kTGV_c(\dom,\mathrm{Sym}^\kTGV(\mathbb{R}^2)), \|\ddiv^\lTGV\, \dualt(x)\|_2 \leq \regup_\lTGV,\,\, \forall x\in\dom \,\,\mbox{and}\,\, \lTGV=0,...,\kTGV-1 \right\}.\notag\end{align}
\end{thm}

Clearly $\DTGV^1_1(u)=\DTV(u)$.

\section{Properties of DTGV} \label{sec:properties}

In this section, we will derive some properties of $\DTGV^{\kTGV}_{\regup}$ on the space of \textit{Bounded Generalized Variation} of order $\kTGV$ ($\BGV^\kTGV$), which is defined as  
\begin{equation*} 
\BGV^\kTGV(\dom) = \left\{\primf\in L^1(\dom) \left| \TGV_\regup^\kTGV(\primf) < \infty \right. \right\}.
\end{equation*}
When $\BGV^\kTGV(\dom)$ is equipped with the norm
\begin{align*}
\|\primf\|_{\BGV^\kTGV(\Omega)} = \|\primf\|_{L^1(\Omega)} + \TGV_\regup^\kTGV(\primf),
\end{align*}
it is a Banach space \cite{Bredies2010}. In the following two propositions we will show that by replacing  $\TGV_\regup^\kTGV$ with $\DTGV_\lambda^\kTGV$ we get an equivalent norm on $\BGV^\kTGV(\dom)$.

\begin{prop}
$\DTGV_\regup^\kTGV: \BGV^\kTGV(\Omega) \to \mathbb{R}$ is a semi-norm.
\end{prop}
\textit{Proof:} Based on the definition of $\DTGV_\regup^\kTGV$, it is obvious that $\DTGV_\regup^\kTGV(u)\geq 0$. Further, with the linearity of the integral we have
\begin{align*}
\DTGV_\regup^\kTGV(\scav\primf) = |\scav|\DTGV_\regup^\kTGV(\primf).
\end{align*}
Define the feasible set of the supremum problem in \cref{def:DTGV} as:
\begin{equation*}
\kset_{D} = \left\{ \dualt\in C^\kTGV_c(\dom,\mathrm{Sym}^\kTGV(\mathbb{R}^2)) \left| \|\ddiv^\lTGV\, \dualt(x)\|_2 \leq \regup_\lTGV,\,\, \forall x\in\dom \,\,\mbox{and}\,\, \lTGV=0,...,\kTGV-1 \right. \right\}.
\end{equation*}
Then, we have that for $\primf_1,\primf_2 \in \BGV_\regup^\kTGV$ 
\begin{align*}
\DTGV_\regup^\kTGV(\primf_1 + \primf_2) &= \sup_{\dualt\in\kset_D} \int_\dom (\primf_1+\primf_2)\ddiv^\kTGV\dualt\dd\scoord \\
& \leq \sup_{\dualt\in\kset_D}  \int_\dom \primf_1\ddiv^\kTGV\dualt\dd\scoord  + \sup_{\dualt\in\kset_D}  \int_\dom \primf_2\ddiv^\kTGV\dualt\dd\scoord \\
& = \DTGV_\regup^\kTGV(\primf_1) + \DTGV_\regup^\kTGV(\primf_2).
\end{align*}
In addition, according to the definition of $\DTGV^{\kTGV}_{\regup}$ we have $\DTGV_\regup^\kTGV(u) = 0$ for any constant function $\primf$. Therefore, we conclude that $\DTGV_\regup^\kTGV$ is a semi-norm on $\BGV^\kTGV(\dom)$. 
\qeda

\begin{prop} \label{prop:equal}
For a function $\primf\in L^1(\dom)$, we have
\begin{align} \label{eqn:eqvrel}
\widE^\kTGV\frac{\min_{\lTGV\in\{0,\kTGV-1\}} \regup_\lTGV}{\max_{\lTGV\in\{0,\kTGV-1\}} \regup_\lTGV} \TGV_\regup^\kTGV(\primf) \leq \DTGV_\regup^\kTGV(\primf) \leq \TGV_\regup^\kTGV(\primf).
\end{align}
\end{prop} 
\textit{Proof:} 
Define the feasible sets of the supremum problems in the definitions of DTGV in \cref{def:DTGV2} and TGV in \cref{eqn:TGV} respectively as
\begin{align*}
\kset_{\el} &= \left\{ \dualt\in C^\kTGV_c(\dom,\mathrm{Sym}^\kTGV(\mathbb{R}^2)) \left| \mathrm{div}^\lTGV\,\dualt \in \regup_l\el^{\widE,\angE}_{\footnotesize\underbrace{2\times \ldots\times 2}_{\kTGV-\lTGV\text{ times}} }(0), \lTGV=0,\ldots,\kTGV-1 \right. \right\}.\\
\kset_{B} &= \left\{ \dualt\in C^\kTGV_c(\dom,\mathrm{Sym}^\kTGV(\mathbb{R}^2)) \left| \mathrm{div}^\lTGV\,\dualt \in \regup_l B_{\footnotesize\underbrace{2\times \ldots\times 2}_{\kTGV-\lTGV\text{ times}} }(0), \lTGV=0,\ldots,\kTGV-1 \right. \right\}.
\end{align*}
Since $a\leq 1$ implies $\el^{\widE,\angE}_{2\times \ldots\times 2}(0) \subseteq B_{2\times \ldots\times 2}(0)$ we see that  $\kset_\el\subset\kset_B$ and hence by the definitions that  $\DTGV_\regup^\kTGV(\primf) \leq \TGV_\regup^\kTGV(\primf)$.

If we shrink the set $\kset_{B}$ to 
\begin{equation*}
\kset_{\widetilde{B}} = \left\{ \dualt\in C^\kTGV_c(\dom,\mathrm{Sym}^\kTGV(\mathbb{R}^2)) \left| \mathrm{div}^\lTGV\,\dualt \in \widE^{\kTGV}\regup_lB_{\footnotesize\underbrace{2\times \ldots\times 2}_{\kTGV-\lTGV\text{ times}} }(0), \lTGV=0,\ldots,\kTGV-1 \right. \right\} ,
\end{equation*}
then we have $\kset_{\widetilde{B}}\subset\kset_\el$. Further, we obtain the inequality
\begin{equation*}
\sup_{ \dualt \in \kset_{\widetilde{B}}} \int_\dom \primf \mathrm{div}^\kTGV\,\dualt \dd \scoord \leq \sup_{ \dualt \in \kset_\el} \int_\dom \primf \mathrm{div}^\kTGV\,\dualt \dd \scoord,
\end{equation*}
that is, $\TGV_{\widE^{\kTGV}\regup}^\kTGV(\primf) \leq \DTGV_\regup^\kTGV(\primf)$. Based on the third statement in proposition 3.3 in \cite{Bredies2010}, we have the relation between TGV-functionals with different weights as:
\begin{equation*}
c \TGV_{\regup}^\kTGV(\primf) \leq \TGV_{\widE^{\kTGV}\regup}^\kTGV(\primf) \quad \mathrm{with}\quad c = \widE^{\kTGV}\frac{\min_{\lTGV\in\{0,\kTGV-1\}} \regup_\lTGV}{\max_{\lTGV\in\{0,\kTGV-1\}} \regup_\lTGV}. 
\end{equation*} 
Hence, $c \TGV_{\regup}^\kTGV(\primf) \leq \DTGV_\regup^\kTGV(\primf)$.
\qeda

There are two straightforward consequences from \cref{prop:equal}. First of all $\BGV^\kTGV(\dom)$ can equivalently be equipped with the norm
\begin{align*}
\|\primf\|_{\BGV_\regup^\kTGV(\Omega)} = \|\primf\|_{L^1(\Omega)} + \DTGV_\regup^\kTGV(\primf).
\end{align*}
Second,  the kernel of $\DTGV^{\kTGV}_{\regup}$ can be characterized:
\begin{cor} \label{cor:DTGVpol}
$\DTGV_\regup^\kTGV(\primf)=0$ if and only if $\primf$ is a polynomial of degree less than $\kTGV$.
\end{cor}
\textit{Proof:} The kernel of $\TGV_\regup^\kTGV=0$ consists of polynomials of degree less than $\kTGV,$ see \cite{Bredies2010}; hence the conclusion follows from \cref{eqn:eqvrel}.  \qeda

In the next proposition we will derive further properties of $\DTGV^{\kTGV}_{\regup}.$


\begin{prop} \label{prop:DTGVh}
$\DTGV_\regup^\kTGV: \BGV^\kTGV \to \mathbb{R}^+\cup\{0\}$ is convex and lower semi-continuous.
\end{prop}
\textit{Proof:} For $\primf_1,\primf_2 \in L^1(\dom)$ and $\scav\in [0,1]$ we have 
\begin{align*}
\DTGV_\regup^\kTGV(\scav\primf_1 + (1-\scav)\primf_2) &= \sup_{\dualt\in \kset_\el} \int_\dom (\scav\primf_1 + (1-\scav)\primf_2)\ \ddiv^\kTGV \dualt\ \dd \scoord \\
&\leq \scav\sup_{\dualt\in \kset_\el} \int_\dom \primf_1\ \ddiv^\kTGV \dualt\ \dd \scoord + (1-\scav)\sup_{\dualt\in \kset_\el} \int_\dom \primf_2\ \ddiv^\kTGV \dualt\ \dd \\
&= \scav\DTGV_\regup^\kTGV(\primf_1) + (1-\scav)\DTGV_\regup^\kTGV(\primf_2).
\end{align*}
Hence, $\DTGV^{\kTGV}_{\regup}$ is convex.

By use of Fatou's lemma we can show the lower semi-continuity of $\DTGV^{\kTGV}_{\regup}$. Let $\{\primf_\intt\}_{\intt\in \mathbb{N}}$ be a Cauchy sequence in $\BGV^{\kTGV}(\dom)$ such that $\primf_{n}\rightarrow\primf$ in $L^1(\dom)$. Based on the definition of $\DTGV^{\kTGV}_{\regup}$ in \cref{def:DTGV2} and Fatou's Lemma, for any $\dualt\in \kset_\el$ we have
\begin{align*}
\liminf_{\intt\to\infty} \DTGV_\regup^\kTGV(\primf_\intt) 
\geq \liminf_{\intt\to\infty} \int_\dom \primf_\intt \mathrm{div}^\kTGV\,\dualt  \dd\scoord
\geq \int_\dom \liminf_{\intt\to\infty} \primf_\intt \mathrm{div}^\kTGV\,\dualt \dd\scoord
= \int_\dom \primf\ \mathrm{div}^\kTGV\,\dualt \dd\scoord.
\end{align*}
Taking the supremum over all $\dualt$ in $\kset_\el$ thus yields
\begin{align*}
\DTGV_\regup^\kTGV(u) \leq \liminf_{\intt\to\infty} \DTGV_\regup^\kTGV(\primf_\intt),
\end{align*}
which means that $\DTGV^{\kTGV}_{\regup}$ is indeed lower semi-continuous. \qeda 
 
In the end of this section we give another similar equivalent definition as in Theorem 3.1 in \cite{Bredies2011} but for $\DTGV_\regup^2$, which is highly attractive in the  numerical implementation.

\begin{thm} \label{thm:DTGV2c}
For 
$\primf \in L^1(\dom)$ we have 
\begin{equation*}
\DTGV^2_\regup(\primf) = \min_{\dualf \in \BD(\dom)} \regup_1 \|\widetilde{\dide} \primf-\dualf \|_{\mathcal{M}} + \regup_0 \| \widetilde{\mathcal{E}} \dualf \|_{\mathcal{M}},
\end{equation*}
where $\BD(\dom)$ denotes the space of vector fields of Bounded Deformation \cite{Temam1985}, $\|\widetilde{\dide}\primf\|_{\mathcal{M}}=\int_\dom \mathrm{d}|\widetilde{\dide}u|=\DTV(u)$, the directional symmetrized derivative $\widetilde{\mathcal{E}}$ 
is the adjoint operator of $\widetilde{div}$ for a vector field.

\end{thm}
Since the proof of \cref{thm:DTGV2c} is following the same lines as the proof of Theorem 3.1 in \cite{Bredies2011} except the change on the divergence operator, for more details we refer the readers to this paper. In addition, in \cref{sec:discnot} we give definitions of all operators in the discrete case.

\section{L$^2$-DTGV$^2_\regup$ Model} \label{sec:rergularization}

In this section we consider \cref{eqn:general} with the regularization term given by $\DTGV_\regup^2,$ i.e.\ 
\begin{equation} \label{eqn:L2DTGV}
\min\limits_{\primf\in \BGV^2(\dom)} \quad\efu(\primf)
\end{equation}
with
\begin{align*}
\efu(\primf) =      \frac{1}{2} \|\fo \primf - \noisef\|_{L^2(\dom)}^2  + \, \DTGV_\regup^2(\primf).
\end{align*}
We call this the  L$^2$-DTGV$^2_\regup$ model. Based on the properties of $\DTGV^2_\regup$ and $\BGV^2(\dom)$, we prove the existence and uniqueness of a solution to \cref{eqn:L2DTGV}.

\begin{thm}
Suppose that $\noisef$ is in $L^2(\dom)$ and $\fo\in\mathcal{L}(L^2(\dom))$ is  injective on the space of affine functions $\affP(\dom)$. Then the $L^2$-DTGV$_\regup^2$ model defined in \cref{eqn:L2DTGV} has a solution. Moreover, the solution is unique.
\end{thm}
\textit{Proof:} Since $\efu$ is bounded from below, we can choose a minimizing sequence $\{\primf_\intt\}_{\intt\in \mathbb{N}}\subset\BGV^2(\dom)$ for \cref{eqn:L2DTGV}. Thus both $\{\|\fo \primf_\intt-\noisef \|_{L^2(\Omega)}\}$ and $\{\DTGV^2_\regup(\primf_\intt)\}$ with $\intt=1, 2, \cdots$ are bounded.

Let $\projP: L^2(\dom) \to \affP(\dom)$ be a linear projection onto the space of affine functions on $\Omega,$ $\affP(\dom)$. Based on the result in proposition 4.1 in \cite{Bredies2011}, we can find a constant $\con>0$ such that
\[
\|\primf\|_{L^2(\Omega)}\leq \con\ \TGV_\regup^2(\primf)
\]
for any $\primf$ in $\text{ker} \projP \subset L^2(\dom)$. Then from \cref{eqn:eqvrel} we have
\begin{align*}
\|\primf\|_{L^2(\Omega)} \leq \widetilde{\con} \ \DTGV_\regup^2(\primf) \qquad \forall \primf \in \text{ker} \projP \subset L^2(\dom),
\end{align*}
with $\widetilde{\con}=\frac{\con\max\{\regup_0, \regup_1\}}{a^2\min\{\regup_0, \regup_1\}}$. By using \cref{cor:DTGVpol} and the triangle inequality on the semi-norm $\DTGV_\regup^2$ we obtain $\DTGV_\regup^2(\primf) = \DTGV_\regup^2(\primf-\projP\primf)$. Hence, we have that $\{\primf_\intt-H\primf_\intt \}_{\intt\in \mathbb{N}}$ is bounded in $L^2(\dom)$. 

Since $\fo$ is injective on the finite-dimensional space $\affP(\dom)$, there is a $\con_1>0$ such that $\|\projP\primf\|_{L^2(\Omega)} \leq \con_1 \|\fo\projP\primf\|_{L^2(\Omega)}$. Further,
\begin{align*}
\|\projP\primf_\intt\|_{L^2(\Omega)} &\leq \con_1 \|\fo\projP\primf_\intt\|_{L^2(\Omega)} \\
&\leq \con_1\left(\|\fo\primf_\intt - \noisef\|_{L^2(\Omega)} + \|\fo(\primf_\intt - \projP\primf_\intt) - \noisef \|_2 \right) \leq \con_2,
\end{align*}
for some $\con_2>0$. It implies that $\{\primf_\intt\}$ bounded in $L^2(\dom).$

Therefore, there exists a subsequence of $\{u_{n}\}$ that converges weakly to a $\primf^*\in L^2(\dom)$. Based on the lower semi-continuity and convexity of $\DTGV_\regup^2$ stated in \cref{prop:DTGVh} 
 we obtain that $\primf^*$ is a minimizer of $\efu$ and hence a solution of the model \cref{eqn:L2DTGV}. 

Based on \cref{prop:DTGVh}, the functional $\efu$ is convex. Furthermore, combining with the result in \cref{cor:DTGVpol} and $\fo$ is injective on $\affP(\dom)$, \cref{eqn:L2DTGV} is strictly convex, thus its minimizer has to be unique \cite{Attouch2006}.
\qeda

\section{Algorithms} \label{sec:algorithms}

In this section we will introduce the algorithms needed for the following numerical simulations. First we will give the notations and discretization of the different operators that our algorithms will require. To keep the notation simple, we re-use the same symbols as in continuous case in the previous sections for the discrete case. Then, we will propose a method for estimating the main direction in images. In the end of the section, we will propose a primal-dual based algorithm for solving the minimization problem in the L$^2$-DTGV$^2_\regup$ model. 

\subsection{Notation and discretization} \label{sec:discnot}
The domain $\dom$ is discretized as an $\pdim$-by-$\pdim$ equidistant pixel-grid with pixel-size $1\times1$. We use $(\xind,\yind)$ to denote a pixel index with $1\leq\xind, \yind\leq M$, such that $u_{\xind,\yind}$ gives the pixel value at $(\xind,\yind)$. Here, for the sake of simplicity we stick to a square pixel-grid, but all proposed algorithms can be easily generalized to any rectangular discretization. 

For $\primf\in \mathbb{R}^{\pdim\times\pdim}$, the discrete gradient operator $\nabla:\mathbb{R}^{\pdim\times\pdim}\to\mathbb{R}^{2\pdim\times\pdim}$ is defined as 
\begin{align*}
\fdm\primf=\begin{pmatrix} \fdm^\forw_{\scoord_1}\primf \\ \fdm^\forw_{\scoord_2}\primf \end{pmatrix},
\end{align*}
where $\nabla^\forw_{\scoord_1}$ and $\nabla^\forw_{\scoord_2}$ are obtained by applying a forward finite difference scheme with symmetric boundary condition, i.e., 
\begin{align*}
(\fdm^\forw_{\scoord_1}\primf)_{\xind,\yind} = \left\{ \begin{array}{ll} \primf_{\xind+1,\yind}-\primf_{\xind,\yind}, & \text{if}\; \xind<\pdim, \\ 0, & \text{if}\; \xind=\pdim, \end{array} \right. \quad \text{and} \quad
(\fdm^\forw_{\scoord_2}\primf)_{\xind,\yind} = \left\{ \begin{array}{ll} \primf_{\xind,\yind+1}-\primf_{\xind,\yind}, & \text{if}\; \yind<\pdim, \\ 0, & \text{if}\; \yind=\pdim. \end{array} \right.
\end{align*}
The divergence operator is defined as the adjoint operator of $\nabla$, i.e., we have $\mathrm{div}=-\nabla^*=(\fdm^\bacw_{\scoord_1}, \fdm^\bacw_{\scoord_2})$, where $\fdm^\bacw_{\scoord_1}$ and $\fdm^\bacw_{\scoord_2}$ utilize the backward finite difference scheme. 

Moreover, based on the relation in \cref{ddiv1d}, the directional divergence for a tensor $\mathbf{v}$ with $\mathbf{v}_{\xind,\yind}=(v^{1}_{\xind,\yind}, v^{2}_{\xind,\yind})^\top$ and $1\leq \xind, \yind\leq \pdim$ can be obtained by calculating 
\begin{align*}
(\ddiv \mathbf{v})_{\xind,\yind} = (\mathrm{div}\ \widetilde{\mathbf{v}})_{\xind,\yind} \quad\mbox{with}\quad \widetilde{\mathbf{v}}_{\xind,\yind} = \rot_{\angE}\scm_{\widE}\mathbf{v}_{\xind,\yind}.
\end{align*}
The corresponding directional gradient operator is $(\widetilde{\grad} \primf)_{\xind,\yind}=\scm_{\widE}\rot_{-\angE}(\grad\primf)_{\xind,\yind}$. 

For a tensor $\dualtm$ with $\dualtm_{\xind,\yind}=\left(\begin{smallmatrix} \dualte_{\xind,\yind}^{11} & \dualte_{\xind,\yind}^{12} \\ \dualte^{12}_{\xind,\yind} & \dualte^{22}_{\xind,\yind} \end{smallmatrix}\right)$ and $1\leq \xind, \yind\leq \pdim$, its divergence can be expressed as 
\[
(\mathrm{div}\, \dualtm)_{\xind,\yind} =  
\begin{pmatrix} (\fdm^\forw_{\scoord_1} \dualte^{11} )_{\xind,\yind} + (\fdm^\forw_{\scoord_2} \dualte^{12} )_{\xind,\yind} \\ (\fdm^\forw_{\scoord_1} \dualte^{12} )_{\xind,\yind} + (\fdm^\forw_{\scoord_2} \dualte^{22} )_{\xind,\yind} \end{pmatrix}. 
\]
Based on the definition in \cref{eqn:ddiv}, the directional divergence, $\ddiv$, can be obtained by calculating
\[
(\ddiv W)_{\xind,\yind}=(\mathrm{div} \widetilde{W})_{\xind,\yind} \quad\mbox{with}\quad \widetilde{W}_{\xind,\yind}=\rot_\angE\scm_\widE W_{\xind,\yind}.
\]

In addition, the directional symmetrized derivative of the tensor $\mathbf{v}$ defined above is given by
\begin{align*}
(\widetilde{\mathcal{E}}\mathbf{v})_{\xind,\yind} = \frac12\left[\scm_{\widE}\rot_{\angE}
\begin{pmatrix} (\fdm^\bacw_{\scoord_1}\dualf^1)_{\xind,\yind} &
( \fdm^\bacw_{\scoord_1}\dualf^2)_{\xind,\yind}\\
(\fdm^\bacw_{\scoord_2}\dualf^1)_{\xind,\yind}  &
(\fdm^\bacw_{\scoord_2}\dualf^2)_{\xind,\yind} \end{pmatrix}+
\begin{pmatrix} (\fdm^\bacw_{\scoord_1}\dualf^1)_{\xind,\yind} &
( \fdm^\bacw_{\scoord_2}\dualf^1)_{\xind,\yind}\\
(\fdm^\bacw_{\scoord_1}\dualf^2)_{\xind,\yind}  &
(\fdm^\bacw_{\scoord_2}\dualf^2)_{\xind,\yind} \end{pmatrix}\rot_{\angE}^\top\scm_{\widE}\right]
\end{align*}
Note that for tensors it still holds that $\widetilde{\mathrm{div}}=-\widetilde{\mathcal{E}}^*$.

\subsection{Detecting the main direction in images} \label{sec:angdir}

In order to apply DTGV as regularization, the parameters $\widE$ and $\angE$ have to be specified. In
this paper, we focus on the case that the texture in images is mainly along one direction. By estimating this main direction, we will obtain the parameter $\angE$. The parameter $\widE$ somehow shows the confidence on the angle estimation. For image restoration, the main direction can be estimated directly from the degraded images. Some classical methods to estimate angles or directions in images could be used here, e.g., 2D Fourier transform and the arctangent function with two arguments. There are also more advanced techniques for angle estimation. Here we list a few of them: the quadrature filter \cite{Granlund1995}, the boundary tensor \cite{Granlund1995}, and the structure tensor \cite{Weickert1998}. Most of these methods estimate the direction pixel-wisely, but in our case we aim for only one main direction for the whole image. In this section, based on our single direction assumption, we will introduce another method for estimating the main direction.

Our direction estimator is inspired by \cite{Setzer2008}, and is mainly composed of three steps. First, we smooth the degraded image in order to reduce the influence of noise. Then, a pixel-wise angle estimation is calculated as 
\begin{align} \label{eqn:PiAEst}
\Theta_{\xind,\yind} =  \left\{ \begin{array}{ll} 0, & \text{if } \eun{(\grad f^{\stdG})_{\xind,\yind}} < 10^{-3}, \\
\arccos \left( \frac{(\grad^\forw_{\scoord_1} f^{\stdG})_{\xind,\yind}}{\eun{(\grad f^{\stdG})_{\xind,\yind}}} \right), & (\grad^\forw_{\scoord_2} f^{\stdG})_{\xind,\yind} \geq 0 \text{ and }  \eun{(\grad f^{\stdG})_{\xind,\yind}} < 10^{-3}, \\ 
2\pi - \arccos \left( \frac{(\grad^\forw_{\scoord_1} f^{\stdG})_{\xind,\yind}}{\eun{(\grad f^{\stdG})_{\xind,\yind}}} \right), &\text{ otherwise,}
\end{array} \right.
\end{align}
where $f^{\stdG}$ denotes the smoothed image from the first step. After that, we smooth the estimated angles in order to remove outliers and features due to noise. At the same time, we introduce the new period for the angles. Note that in \cite{Setzer2008} the focus is on restoring rectangular shapes, therefore $\pi/2$-period is used. But in our case, we need a $\pi$-period estimate. Moreover, we only need the main direction in the image, which is obtained by calculating the mean over the pixel-wise angles. In \cref{alg:ThEst} the details of the main direction estimation method are outlined. It should be noted that if we do not calculate the mean of the angles, we will have pixel-wise angle estimates, which can be utilized for restoring images with multiple angles in the future. Since in this paper we mainly focus on the analysis of DTGV under continuous setting and its extension to spatially varying angles is not trivial, the applications on restoring images with multiple angles are outside the scope of this work.

\begin{algorithm}[htbp]
\caption{Main Direction Estimation Method}
\label{alg:ThEst}
\begin{algorithmic}[1]
\State Input smoothing parameter $\smp$ and the degraded image $f$.
\State Smooth the degraded image by implementing Gaussian blur: $f^\stdG=\kerG f$, where $\kerG$ denotes Gaussian blurring operator with mean 0 and variance $\sigma^{2}$.
\State Estimate pixel-wise direction $\Theta^\aone$ according to \cref{eqn:PiAEst}.
\State Introduce $\frac{\pi}{2}$-period to the angles and smooth them:\\[2mm]
$\cosm_{i,j}^\aone=\cos(4\Theta_{i,j}^\aone),$\\ 
$\sinm_{i,j}^\aone=\sin(4\Theta_{i,j}^\aone),$\\
$(\cosm^\atwo,\sinm^\atwo) = \displaystyle{\arg\min_{\cosm^\atwo,\sinm^\atwo} \sum_{\xind,\yind} \eun{(\grad f)_{\xind,\yind} }^2 
\eun{\left(\begin{array}{c}\cosm^\aone_{\xind,\yind}  \\ \sinm^\aone_{\xind,\yind} \end{array}\right)-\left(\begin{array}{c}\cosm^\atwo_{\xind,\yind}  \\ \sinm^\atwo_{\xind,\yind} \end{array}\right)}^2 
+ \smp \left( \eun{(\grad \cosm^\atwo)_{\xind,\yind}}^2 + \eun{(\grad \sinm^\atwo)_{\xind,\yind}}^2 \right)}$. \\[2mm]
The smoothed pixel-wise angles, $\Theta^\atwo$, are obtained by implementing \cref{eqn:PiAEst} with $\left(\begin{smallmatrix} \cosm^\atwo\\\sinm^\atwo\end{smallmatrix}\right)$ instead of $\grad f^{\stdG}$ as input.
\State Calculate the main direction: $\angE=\frac{1}{|\dom|} \displaystyle\sum_{\xind,\yind} \Theta^\atwo_{\xind,\yind} $. 
\State Transform into $\pi$-periodic angle according to\\
$\angE \leftarrow
\left\{ \begin{array}{ll} -\angE, & \text{if }  \displaystyle (\grad^{\forw}_{\scoord_1} f^\stdG) : (\grad^{\forw}_{\scoord_2} f^\stdG) \leq 0, \\ 
\frac{\pi}{2}-\angE, & \text{otherwise,} \end{array} \right. $\\
where ``:'' denotes the Frobenius inner product.
\Return $\angE$.
\end{algorithmic}
\end{algorithm}

\subsection{The Chambolle-Pock algorithm} \label{sec:CPDTGV}
Corresponding to \cref{eqn:L2DTGV} we formulate the discrete L$^2$-DTGV$_\regup^2$ model as
\begin{equation} \label{eqn:discL2DTGV} 
\underset{\primf\in\mathbb{R}^{\pdim\times\pdim}}{\min}\ \mathcal{J}(u):= \frac{1}{2}\|\sysm\primf-\noisef\|_F^2 + \DTGV^2_\regup(\primf),
\end{equation}
where $\noisef\in\mathbb{R}^{\pdim\times\pdim}$ and $A: \mathbb{R}^{\pdim\times\pdim}\rightarrow\mathbb{R}^{\pdim\times\pdim}$ denotes the identity operator (denoising problem) or a blurring operator (deblurring problem).
Since the minimization problem in \cref{eqn:discL2DTGV} is convex, many optimization algorithms could be used to solve it, e.g. Nesterovs method \cite{Nesterov1983}, the FISTA algorithm \cite{Beck2009a}, the alternating direction method with multipliers (ADMM) \cite{Boyd2010}, and any of the many primal-dual-based methods. Here, due to the simplicity of the implementation, we utilize the Chambolle-Pock primal-dual (CP) algorithm \cite{Chambolle2011a} to solve our problem. 

Referring to the algorithm proposed in \cite{Bredies2014}, we can rewrite the data-fitting term in \cref{eqn:discL2DTGV}, $\funcB(u) = \frac{1}{2} \|\sysm \primf - \noisef\|_F^2$, as 
\begin{align} \label{eqn:CPreform1}
\funcB(\primf) = \max_{\imq\in\mathcal{U}}\ \langle \sysm \primf,\imq\rangle - \frac{1}{2}\|\imq\|_F^2-\langle\noisef,\imq\rangle,
\end{align}
where $\mathcal{U}=\mathbb{R}^{\pdim\times\pdim}$. Combining with the result in \cref{thm:DTGV2c}, we obtain the primal-dual formulation of \cref{eqn:discL2DTGV}
\begin{align*}
\min_{\primf\in\mathcal{U},\mathbf{v}\in\mathcal{V}}\ \max_{\imq\in\mathcal{U},\mathbf{p}\in\mathcal{P},\dualtm\in\mathcal{W}}\ 
&\langle \sysm\primf,\imq \rangle
- \frac{1}{2} \|\imq \|_F^2 
- \langle \noisef,\imq \rangle
+ \langle\widetilde{\grad} \primf -\mathbf{v}, \mathbf{p} \rangle
+ \langle\widetilde{\mathcal{E}}\mathbf{v},\dualtm \rangle
\end{align*}
where $\mathcal{V}=\mathbb{R}^{2\pdim\times\pdim}$, $\mathcal{P}=\{\mathbf{p}: \dom\rightarrow\mathbb{R}^{2}\ |\ \|\mathbf{p}_{i,j}\|_{2}\leq\regup_{1} \mbox{ for } \forall (i,j)\in\dom\}$, $\mathcal{W}=\{W: \dom\rightarrow \mathrm{Sym}^{2}(\mathbb{R}^{2})\ |\ \|W_{i,j}\|_{F}\leq\regup_{0} \mbox{ for } \forall (i,j)\in\dom\}$. This is a generic saddle-point problem, and we can apply the CP algorithm proposed in \cite{Chambolle2011a} to solve it. The algorithm is summarized in Algorithm 2.

\begin{algorithm}
\caption{The CP algorithm for solving L$^2$-DTGV$^{2}_{\regup}$}
\label{alg:CPDTV}
\begin{algorithmic}[1]
\State Require $\noisef$, $\sysm$, $\regup$, $\widE$, $\angE$ and $\tolr$.
\State Estimate Lipschitz constant $\lipc$, e.g. using power-method for $\sysm$.
\State Initialize $\primf^0=\bar{\primf}^0=0$, $\dualf^0=0$, $\imq^0=0$, $\dualCP^0=\bar{\dualCP}^0=0$, $\dualte^0=0$, $\odiff^0=0$, $\stepA<\frac{1}{\sqrt{\lipc}}$, $\stepB<\frac{1}{\sqrt{\lipc}}$. 
\State Run loop until stopping criterion is met: \\
\textbf{while} $\odiff^\niter > \tolr$ \textbf{do}
\begin{align*}
\mathbf{p}^{\niter+1}&= \arg\max_{\mathbf{p}\in\mathcal{P}} \ \langle\widetilde{\grad} \bar{\primf}^{\niter} -\bar{\mathbf{v}}^{\niter}, \mathbf{p} \rangle-\frac{1}{2\stepA} \|\mathbf{p}-\mathbf{p}^{\niter}\|^2_{F}\\
&=\projS_{\regup_1} \left( \mathbf{p}^\niter + \stepA\left( \widetilde{\grad} \bar{\primf}^\niter - \bar{\mathbf{v}}^\niter \right)  \right)\\
\dualtm^{\niter+1} &=\arg\max_{\dualtm\in\mathcal{W}}\ \langle\widetilde{\mathcal{E}}\bar{\mathbf{v}}^{\niter},\dualtm \rangle-\frac{1}{2\stepA}\|\dualtm-\dualtm^{\niter}\|^2_{F}\\
&= \projS_{\regup_0} \left( \dualtm^\niter + \stepA\widetilde{\mathcal{E}} \bar{\mathbf{v}}^\niter \right)\\
\imq^{\niter+1} &=\arg\max_{\imq\in\mathcal{U}}\ \langle \sysm \bar{\primf}^{\niter},\imq \rangle - \frac{1}{2} \|\imq \|_F^2 - \langle \noisef,\imq \rangle -\frac{1}{2\stepA} \|\imq-\imq^{\niter}\|^2_{F}\\
&= \frac{1}{1+\stepA} \left(\imq^\niter + \stepA (\sysm\bar{\primf}^{\niter}-\noisef) \right) \\
\primf^{\niter+1} &=\arg\min_{\primf\in\mathcal{U}}\ \langle \sysm\primf,\imq^{\niter+1} \rangle+ \langle\widetilde{\grad} \primf, \mathbf{p}^{\niter+1} \rangle +\frac{1}{2\stepB} \|\primf-\primf^{\niter}\|^2_{F}\\
&= \primf^\niter + \stepB \left(\ddiv \mathbf{p}^{\niter+1} - \sysm^* \imq^{\niter+1} \right)  \\
\mathbf{v}^{\niter+1}&=\arg\min_{\mathbf{v}\in\mathcal{V}}\  -\langle \mathbf{v}, \mathbf{p}^{\niter+1} \rangle + \langle\widetilde{\mathcal{E}}\mathbf{v},\dualtm^{\niter+1} \rangle +\frac{1}{2\stepB} \|\mathbf{v}-\mathbf{v}^{\niter}\|^{2}_{F}\\
&=  \mathbf{v}^\niter + \stepB \left( \mathbf{p}^{\niter+1} + \ddiv\ \dualtm^{\niter+1} \right) \\
\bar{\primf}^{\niter+1}&=2\primf^{\niter+1} - \primf^\niter\\
\bar{\mathbf{v}}^{\niter+1} &= 2\mathbf{v}^{\niter+1} - \mathbf{v}^\niter \\
\odiff^{\niter+1}&= \frac{\left|\efu\left(\primf^{\niter}\right)-\efu\left(\primf^{\niter+1}\right)\right|}{\efu\left(\primf^{\niter}\right)}
\end{align*}
\textbf{end while} \\
\Return $\primf^{\niter+1}$.
\end{algorithmic}
\end{algorithm}

In \cref{alg:CPDTV}, $\stepA$  and $\stepB$ denote the dual and primal step-sizes, respectively. In addition, $\projS_{\regup}$ is a set-projection operator defined as
\begin{align*}
[\projS_\regup(\xi)]_{i,j} = \frac{\xi_{i,j}}{\max\left(1,\frac{|\xi_{i,j}|}{\regup}\right)}.
\end{align*}
If $\xi\in\mathcal{P}$, then $|\xi_{i,j}|$ is with 2-norm; and if $\xi\in\mathcal{W}$, then $|\xi_{i,j}|$ is with Frobenius norm. Here, we use the relative changes of the objective function in \cref{eqn:discL2DTGV} to define the stopping criterion, since the objective function is essentially what we desire to minimize and it is simple to calculate.

\section{Numerical Experiments} \label{sec:numexp}
In this section, we provide numerical experiments to study the behavior of our method. First, we examine the direction estimation method proposed in \cref{alg:ThEst} for a series of noisy images, ranging from low- to high-level noise. Since the directional regularization requires additional parameters $(\angE,\widE)$, we then empirically examine how these parameters influence the solution of \cref{eqn:discL2DTGV}. In the end, we compare $\DTGV^2_\regup$- and DTV-regularization with $\TGV_\regup^2$- and classical TV-regularization for some denoising and deblurring problems, where the blurring and noise have been simulated. Based on the balance between the computational complexity and restoration improvement, the second order of $\TGV$ is mostly widely used, so in our numerical experiments we only investigate $\DTGV^2_\regup$ and compare it with $\TGV_\regup^2$. To simplify the notations, we refer to $\TGV_\regup^2$ and $\DTGV^2_\regup$ as TGV and DTGV in this section. In addition, same as in \cite{Bredies2010} we fix the ratio between the two regularization parameters, i.e., $\frac{\regup_0}{\regup_1}=2$, which commonly yields good restoration results for images. The tolerance in the \cref{alg:CPDTV} has been chosen as $10^{-6}$ and all simulated experiments are implemented in Matlab R2016a.

\subsection{Robustness of direction-estimation}
Since the estimation of the main direction plays an important role in our method, we first demonstrate the performance of the presented direction-estimation method in \cref{alg:ThEst}. Here, we use one simulated image and one real image for numerical experiments, and test our direction-estimation method on the images with up to 50\% Gaussian noise. In all tests, we set $\sigma=10$ for the smoothing step, i.e., Step 2 in \cref{alg:ThEst}. The numerical results are shown in \cref{fig:test_th_est01}. It is clear that for both test images, until the noise level (nl) reaching to $20\%$, our method provides estimates within $\pm 15^\circ$ of the main direction. Especially for the real image, even with $50\%$ noise, the estimation is still very accurate.

\figone{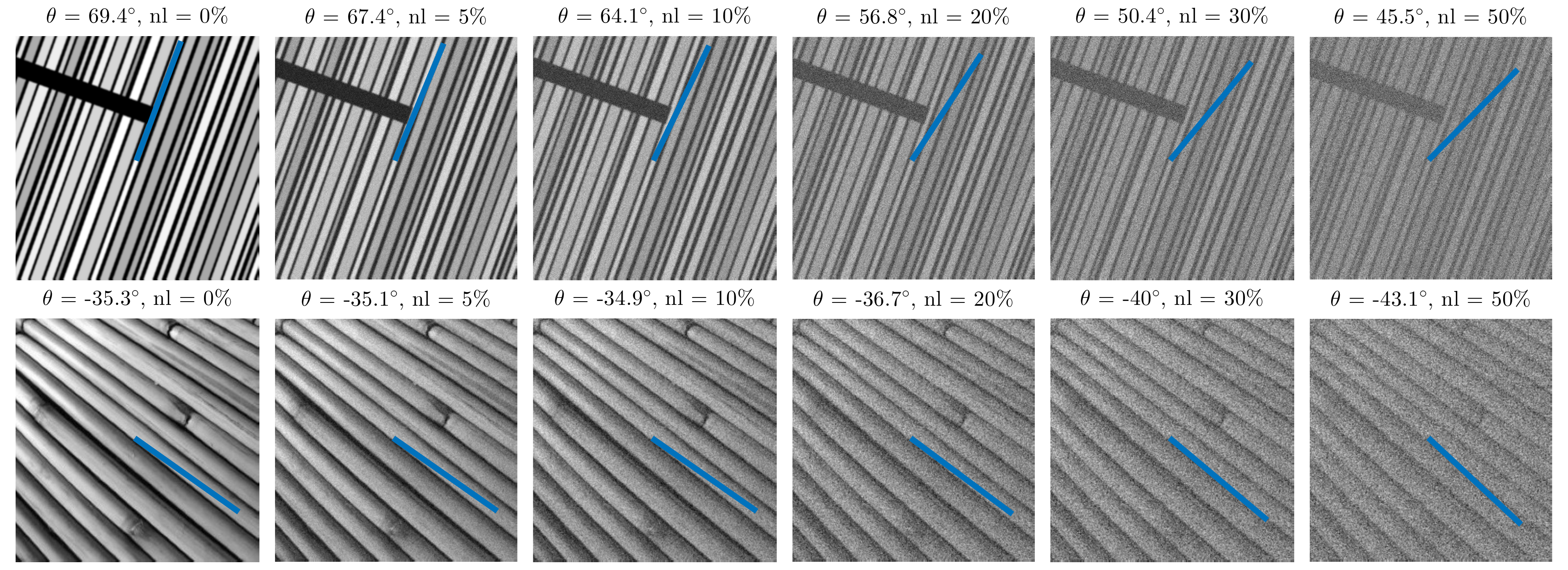}{1}{Estimating main direction in test images with increasing noise-level (nl). The top image size is 512$\times$512 and bottom one is 253$\times$253. Main direction estimate $\angE$ is written in degrees and visualized on top of the noisy images with a blue line from the center.}

\subsection{Role of DTGV parameters}
In the DTGV functional, besides the parameters $\regup_{0,1}$ there are another two important parameters: $\angE$ and $\widE$, where $\angE$ is an angle and $\widE$ determines the ratio of anisotropy. If $\widE=1$, DTGV becomes identical to the rotation invariant TGV. In this section, we will test the influence of these two parameters when DTGV is used as a regularizer, and also seek the robustness of the L$^2-\DTGV$ solutions with respect to the parameter choices.

\figone{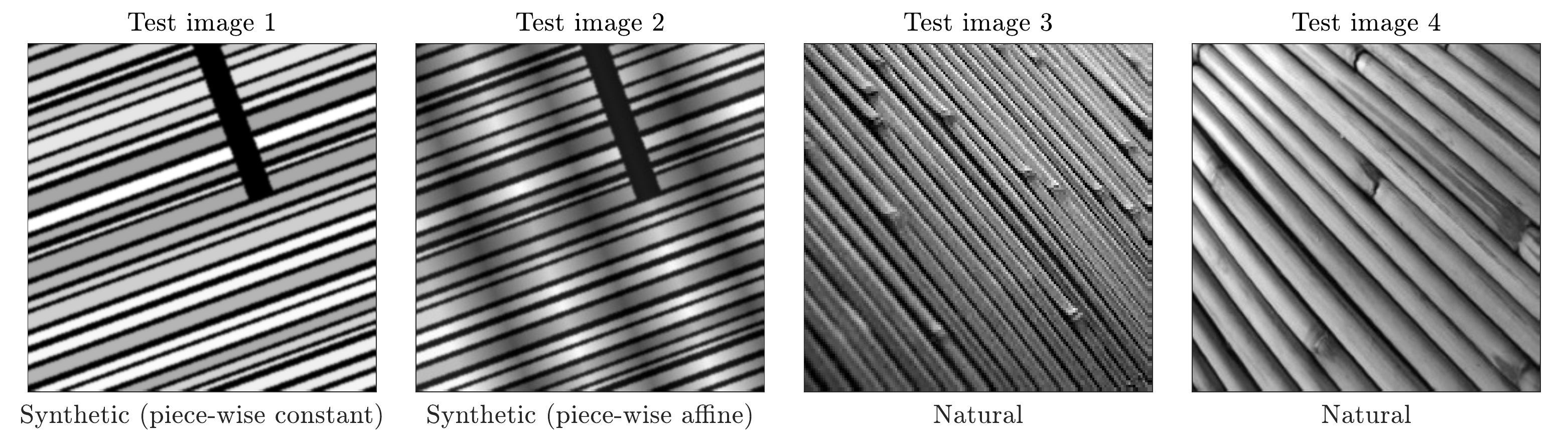}{0.9}{Ground truth images used in the numerical tests.}
\figone{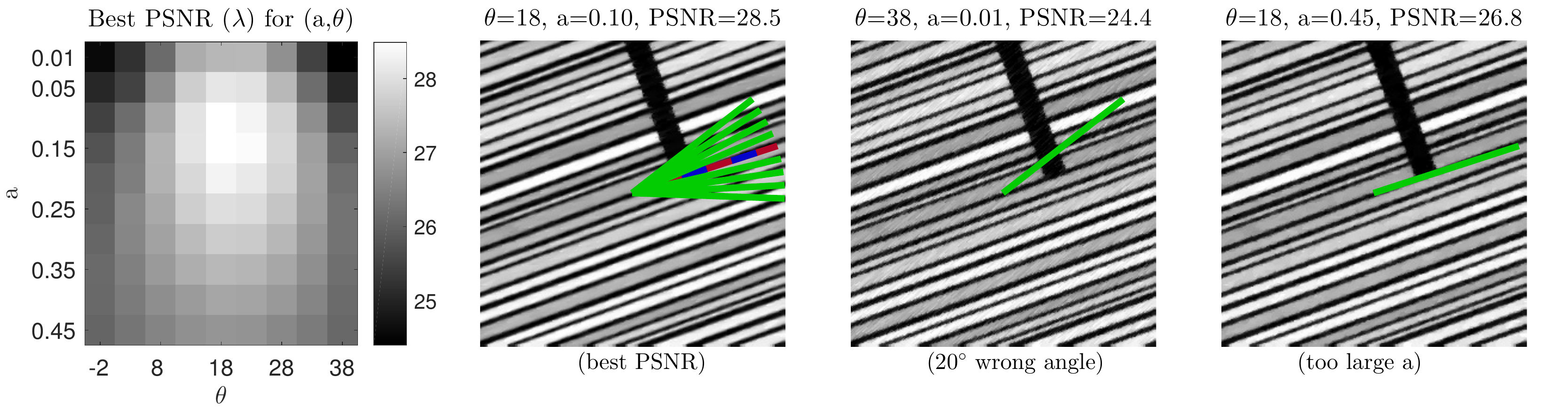}{1}{Test image 1: piece-wise constant image with size 256-by-256. DTV-regularized denoising problem. Lines mark specific angles: {\color{pGreen} tested angles}, {\color{pRed} estimated angle} and {\color{pBlue} best angle}.}
\figone{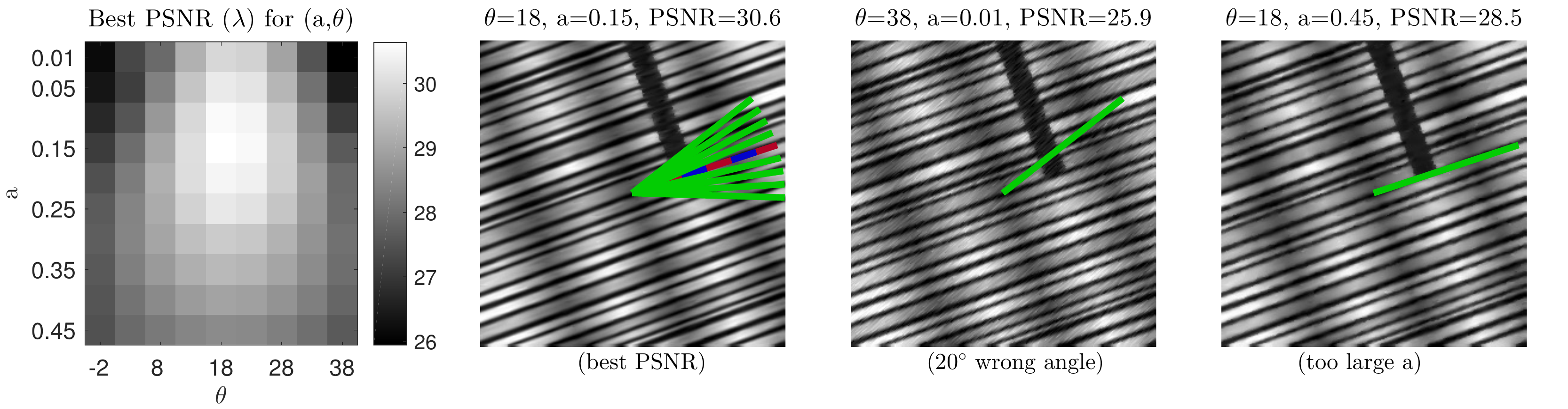}{1}{Test image 2: piece-wise affine image with size 256-by-256. DTGV-regularized denoising problem. Lines mark specific angles: {\color{pGreen} tested angles}, {\color{pRed} estimated angle} and {\color{pBlue} best angle}.}
\figone{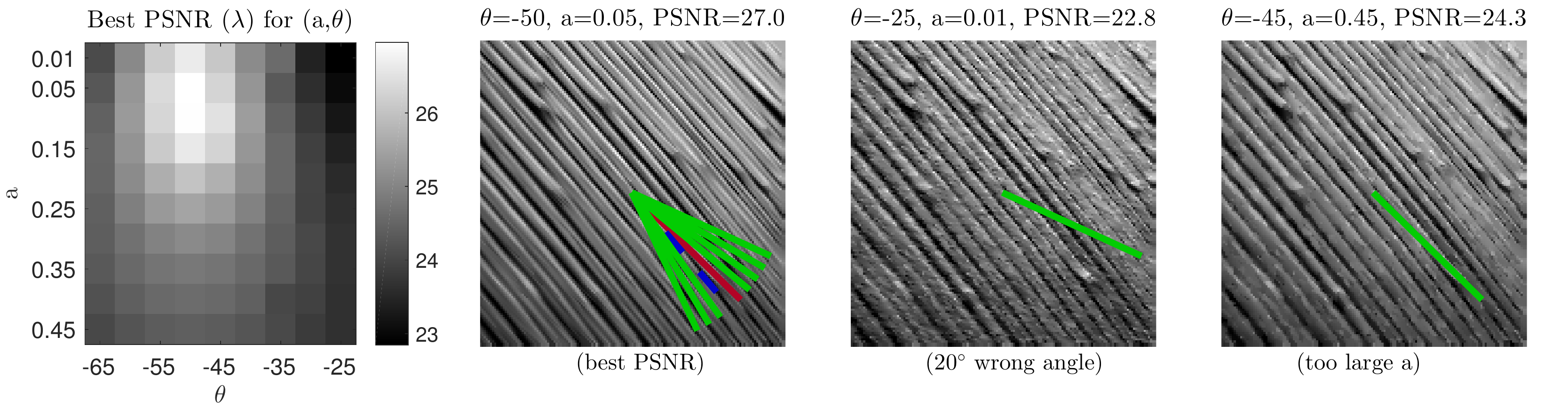}{1}{Test image 3: real image with size 145-by-145. DTGV-regularized denoising problem. Lines mark specific angles: {\color{pGreen} tested angles}, {\color{pRed} estimated angle} and {\color{pBlue} best angle}.}
\figone{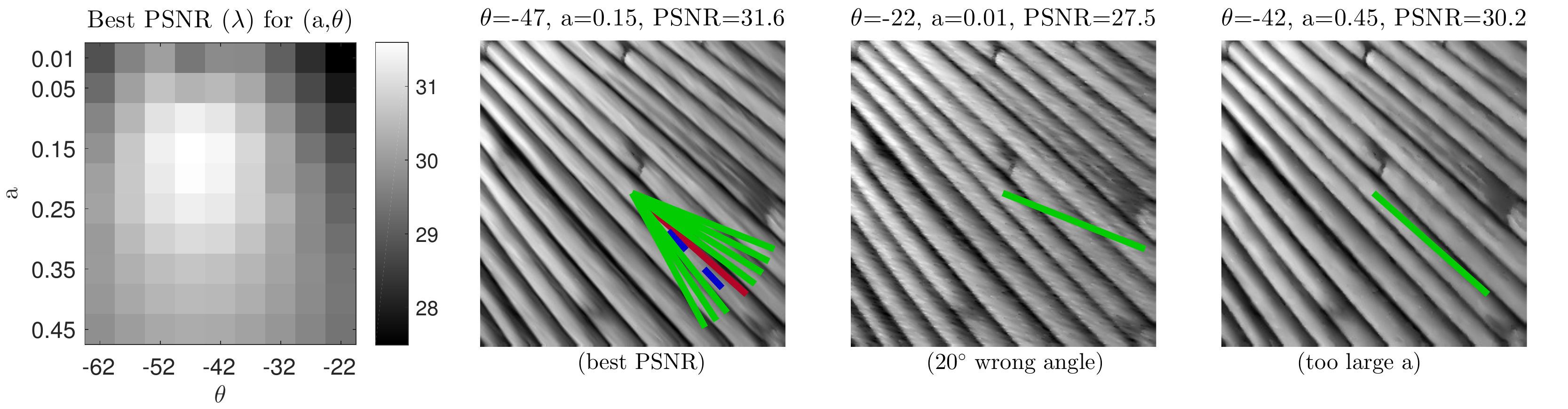}{1}{Test image 4: real image. DTGV-regularized denoising problem. Lines mark specific angles: {\color{pGreen} tested angles}, {\color{pRed} estimated angle} and {\color{pBlue} best angle}.}

In order to study the influence from $\widE$ and $\angE$ solely, in each test with fixed values for $(\widE, \angE)$ we adjust the regularization parameter $\regup$ and pick the one that provides the highest peak signal-to-noise ratio (PSNR) value. The test images used for our numerical algorithm consist of two simulated images and two natural images as shown in \cref{fig:test_ims}. Since the test image in \cref{fig:par_script03_b_DTV} is piece-wise constant, DTV is used as regularizer instead of DTGV. All test images are corrupted by 10\% Gaussian noise, and the operator $\sysm$ is set as the identity. In \cref{fig:par_script03_b_DTV}-\ref{fig:par_script02_b_DTGV}, we visualize the PSNR values for different choice of $(\widE, \angE)$, where we test $\widE\in[0.01,0.45]$ and $\angE\in[\bar{\angE}-20^\circ, \bar{\angE}+20^\circ]$ with $\bar{\angE}$ as the estimated main direction by \cref{alg:ThEst}. In addition, for each test image we also show three restoration results: the best one, according to PSNR, one with $\angE=\bar{\angE}+20^\circ$ and small $\widE$, and one with $\angE=\bar{\angE}$ and large $\widE$.

Obviously, the use of directional regularization improves the PSNR values significantly when $\angE$ coincides with the main direction of the image. Moreover, with a good direction estimation the highest PSNR values are usually achieved by choosing a relatively small $\widE$. From the restored images we can see that with incorrect $\angE$ and a small $\widE$ there are some line artifacts along the direction of $\angE$. The reason is that with a small $\widE$ the textures in the images are forced to be restored along the incorrect direction. With correct $\angE$ and a large $\widE$ the restored results look very similar to the TV- or TGV-regularized results. Especially in the last figure of \cref{fig:par_script03_b_DTV}, staircasing artifacts start to appear.

When DTGV is used as the regularizer to penalize the textures that are not oriented along the main direction of the images, we take $\angE$ as the main direction and choose $\widE$ depending on how much to penalize the textures that are not along $\angE$. Hence the selection of $(\angE, \widE)$ only depends on the orientation of the textures and not on the noise level. In order to empirically confirm the noise robustness, in \cref{fig:par_script02_d_DTGV} we give an example similar to the one in \cref{fig:par_script02_b_DTGV}, but now with the image corrupted by 20\% additive Gaussian noise. In this test, we still adjust the regularization parameter $\regup$ and pick the one that provides the highest PSNR value. Note that $\regup$ controls the trade-off between a good fit to the noisy image and the smoothness from DTGV, so it varies due to the different noise level. We can see that the PSNR figure with respect to $(\widE, \angE)$ in \cref{fig:par_script02_d_DTGV} is very similar to the one in \cref{fig:par_script02_b_DTGV} except the different range of the PSNR values, which demonstrates that the choice of $(\widE, \angE)$ is independent on the noise levels.

\figone{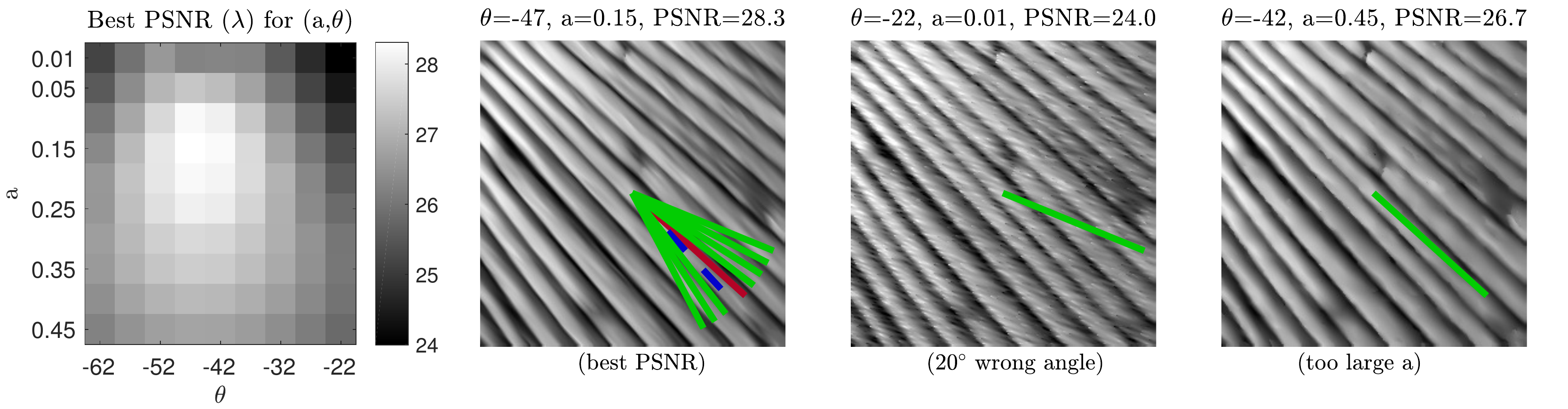}{1}{Test image 4: real image. DTGV-regularized denoising problem (with 20\% Gaussian noise). Lines mark specific angles: {\color{pGreen} tested angles}, {\color{pRed} estimated angle} and {\color{pBlue} best angle}.}

The tests in \cref{fig:par_script03_b_DTV}-\ref{fig:par_script02_d_DTGV} can also serve as robustness tests for our method. They show that our method is robust with respect to the parameters for $\widE\in[0.1, 0.2]$ and $\angE$ within $[\bar{\angE}-5^\circ, \bar{\angE}+5^\circ]$. 
So in the following numerical experiments, we will use \cref{alg:ThEst} to estimate $\angE$ and fix $\widE$ to $0.15$.

\subsection{Image denoising}

To show the improvement of imposing the direction prior into the regularizer we empirically compare DTGV and DTV with TGV and TV for image denoising problems. We use four different test images with 10\% and 20\% Gaussian noise respectively. In all tests, after many experiments with different choices of the regularization parameter $\regup$, the ones that give the best PSNRs are presented here. 

Comparing the results from the four different regularization techniques in \cref{fig:comp_n10} and \cref{fig:comp_n20}, we see that both visually and quantitatively in terms of PSNR the improvement by directional regularizers is evident. Especially the PSNR values of the results from solving the L$^2$-$\DTGV$ model are at least 2.2dB higher than the ones from the L$^2$-$\TGV$  results. The textures in the images are obviously much better preserved by using DTV and DTGV than by using the two isotropic regularizers, i.e., TV and TGV. Note that in both results from TV and DTV for the second test image, the staircasing artifact is observable along the main direction. This is due to the test image being piece-wise affine, while the TV and DTV regularizers are based on an assumption of piece-wise constant images. In this case, by using DTGV as regularizer the artifacts are successfully removed and the textures are well preserved. Hence, these tests demonstrate the advantages of including directional information in the regularizer when dealing with directional images.

In the first two simulated test images, we specially add a dark region perpendicular to the main direction in order to demonstrate the potential artifacts from the directional regularizers. We can see that near the boundary of the dark region some artifacts in the restored results by DTV and DTGV appear. These artifacts are due to the diffusion of the different intensities along the main direction. In addition, the similar artifacts can be observed in the results from DTV and DTGV for restoring the two real images, especially for lines that are close to perpendicular to the main direction.  In order to clearly show these artifacts, in \cref{fig:comp_n20_2_zoom} we zoomed in on the center of the last test image and show the comparison of the four restorations. These artifacts mainly occur when the direction prior is not met, e.g. a part of the image textures do not follow the main direction.

In the $\DTGV^2_\regup$ regularizer, the directional information is introduced both in the first and second-order derivatives. To demonstrate the advantage of this we compare our method with the anisotropic TGV method, ITGV, proposed in \cite{Ranftl2012,Ferstl2013}. In ITGV, a diffusion tensor is estimated from the noisy image and used only in the first-order derivative term. Since in this paper we only consider the case of textures with one global direction, in order to have a fair comparison in \cref{fig:ITGV_comp} we  give the result from ITGV method with the same direction estimation as ours; we mark this result by the name DITGV.  Since in ITGV the diffusion tensor is highly influenced by noise, its performance is limited. Comparing the results from DITGV and DTGV, we can clearly see that the advantage of including the directional information also in the second-order derivative terms.

\figone{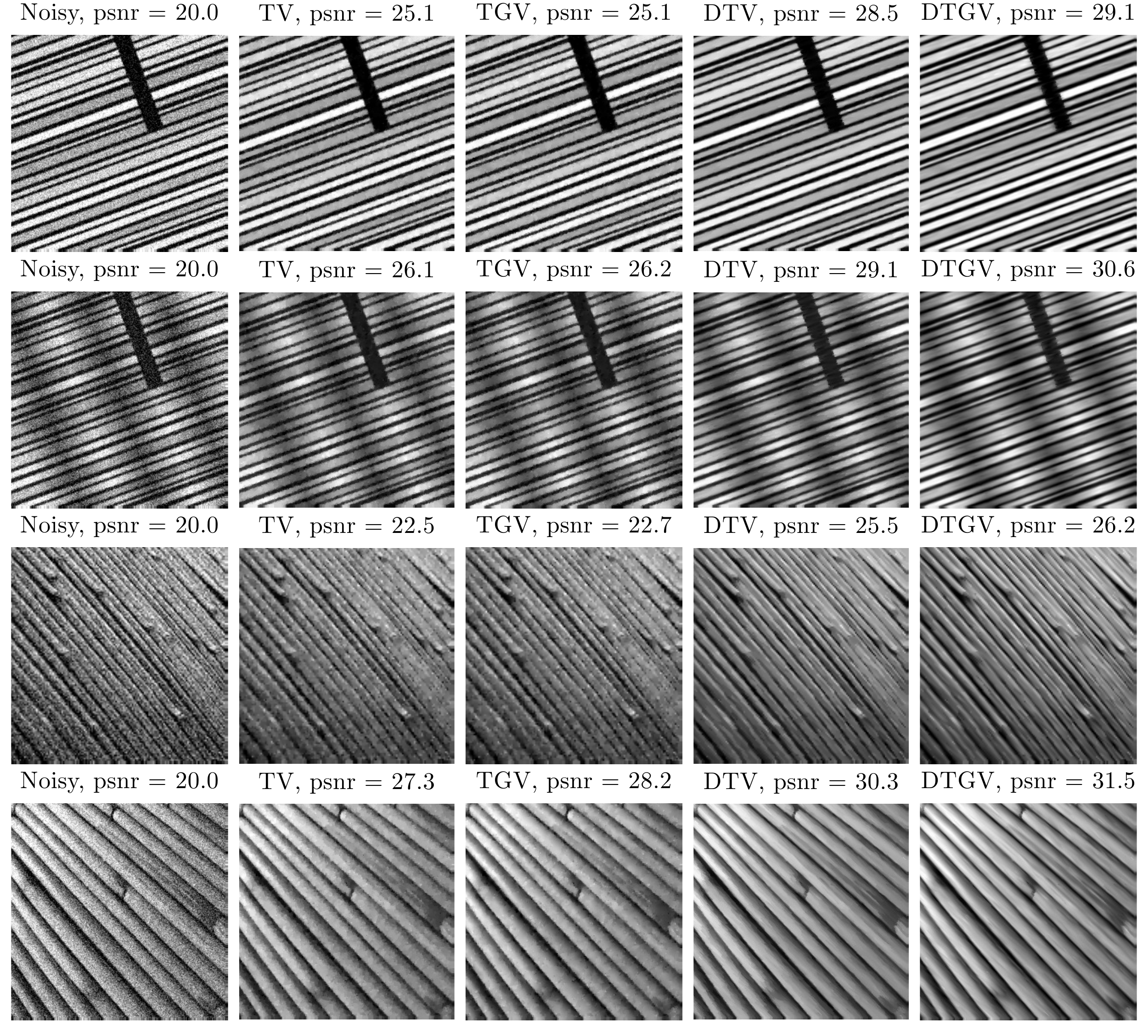}{1}{Comparison of TV, TGV, DTV and DTGV regularizers for four different denoising problems with 10$\%$ Gaussian noise.}
 
\figone{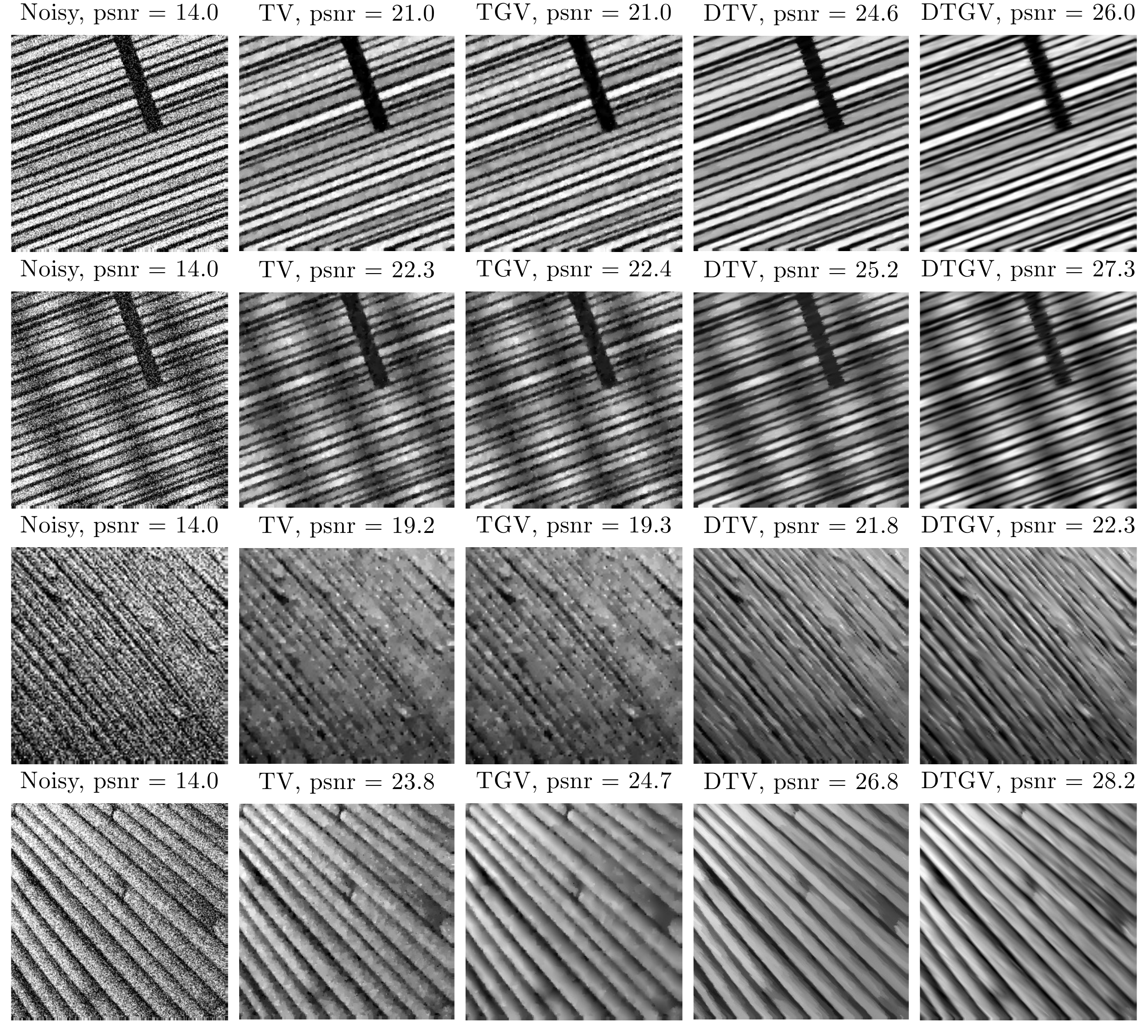}{1}{Comparison of TV, TGV, DTV and DTGV regularizers for four different denoising problems with 20$\%$ Gaussian noise.}

\figone{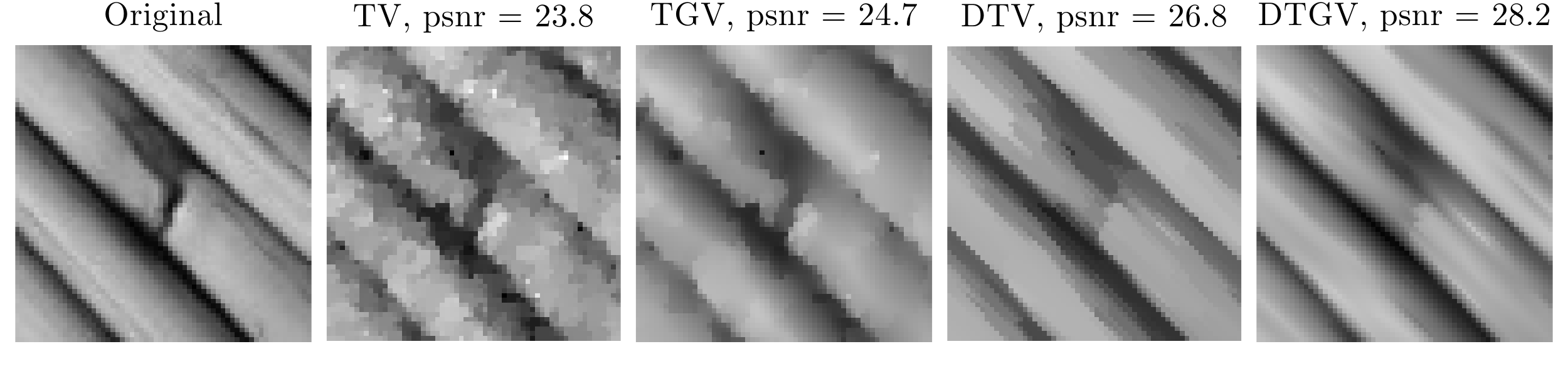}{1}{The zoomed-in regions of the restored results shown in the last row of \cref{fig:comp_n20}.}

\figone{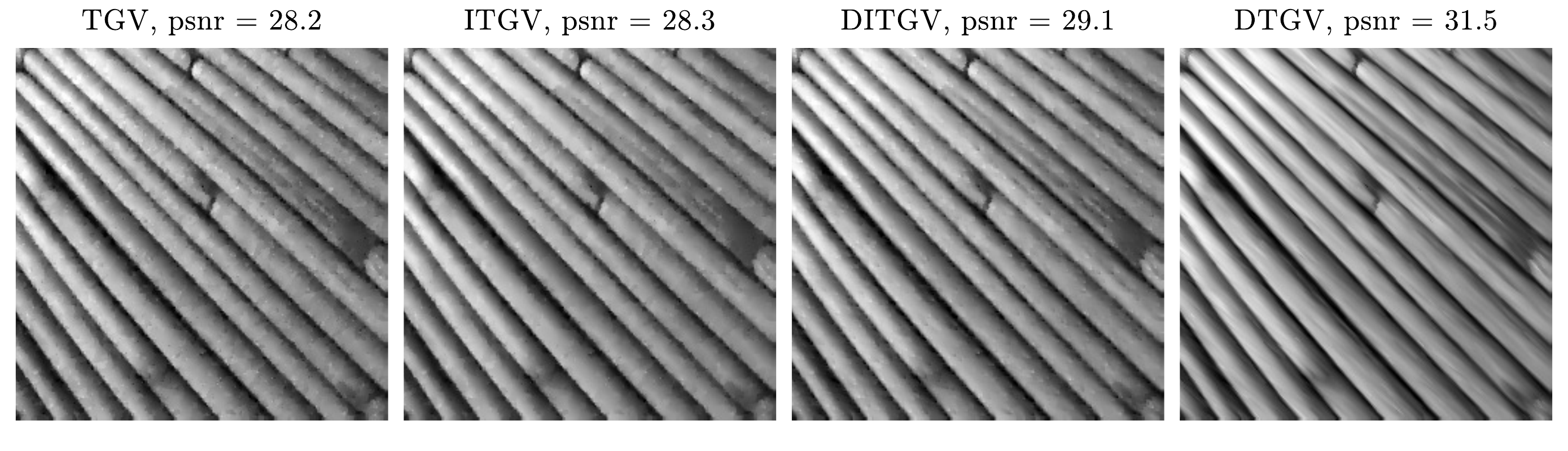}{0.95}{Comparison of TGV, ITGV, DITGV and DTGV regularizers for denoising problem with 10\% Gaussian noise.}

\subsection{Image deblurring and denoising}

In order to test directional regularization on a more complicated problem, we consider the restoration of noisy blurred images. In our experiment, the blurring operator is s et as Gaussian blur with a standard deviation of 2. Further, after being blurred, the test image is corrupted by 10\% Gaussian noise. In our method, we still use \cref{alg:ThEst} to estimate the main direction in the image.

In \cref{fig:comp_blur2_c}, we show the degraded image and the restored results by using TV, TGV, DTV and DTGV regularizers. It is clear that both TV and TGV cannot help to restore the edges correctly, while the methods with the directional regularizers are much better at restoring the textures and removing the noise. In addition, the DTV result is heavily influenced by staircasing artifacts, since this test image is piece-wise affine which does not fit with the piece-wise constant assumption for the method. Due to the use of higher order derivatives, those artifacts do not appear in the DTGV result, as expected. Quantitatively, the PSNR value is increased by at least 2.5dB when the directional regularizer is utilized.  

\figone{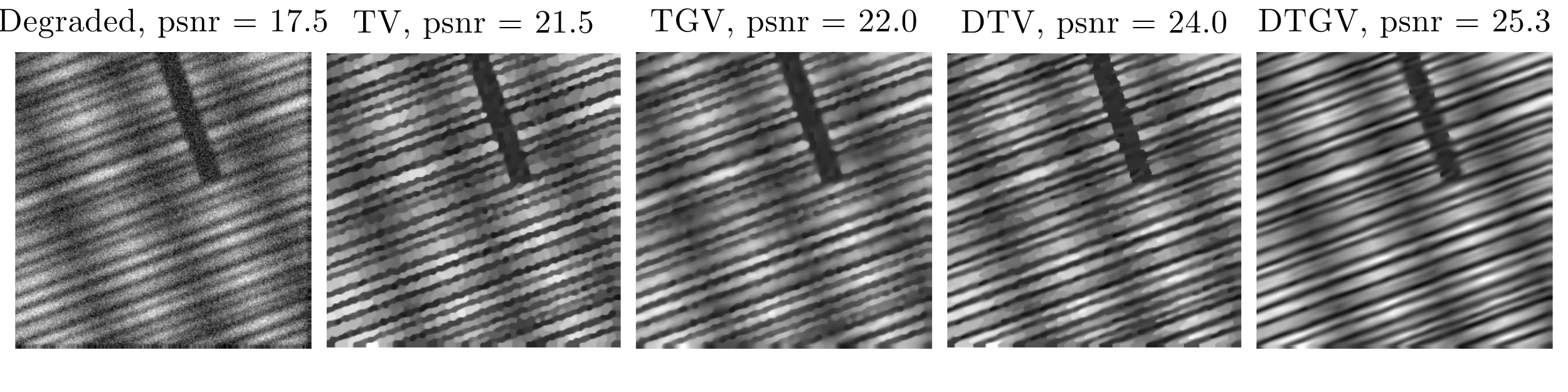}{1}{Comparison of TV, TGV, DTV an DTGV regularizers for a piece-wise affine noisy blurred image.}

\section{Conclusion} \label{sec:conclusion}

In this paper, we propose a new directional functional, directional total generalized variation (DTGV), and study its mathematical properties. Then we combine DTGV with the least-squares data-fitting term and propose a new variational model, L$^2$-$\DTGV$, for restoring images whose textures mainly follow one direction. We prove the existence and uniqueness of a solution to our proposed model, and apply a primal-dual algorithm to solve the minimization problem. Since the new proposed DTGV functional requires the input of the main direction of the images, we also propose a direction estimation algorithm, which can be easily extended to spatially varying direction estimation. Numerical results show that the direction estimation algorithm is reliable and the improvement for restoring directional images by using DTGV as regularizer is significant compared to using isotropic regularizers. In addition, we also try to discover the potential artifacts from DTGV. In order to reduce the artifacts from DTGV, we intend to extend our method to deal with multi-directions or spatially varying directions in the future.

\section*{Acknowledgements}
The work was supported by Advanced Grant 291405 from the European Research Council and grant
no. 4002-00123 from the Danish Council for Independent Research — Natural Sciences.

\bibliographystyle{abbrv}
\bibliography{DTVproject}

\begin{thebibliography}{10}

\bibitem{Attouch2006}
H.~Attouch, G.~Buttazzo, and G.~Michaille.
\newblock {\em {Variational Analysis in Sobolev and BV Spaces: Applications to
  PDEs and Optimization}}.
\newblock Number~1. 2006.

\bibitem{Aubert2006}
G.~Aubert and P.~Kornprobst.
\newblock {\em {Mathematical problems in image processing: partial differential
  equations and the calculus of variations}}.
\newblock 2006.

\bibitem{Bayram2012}
I.~Bayram and M.~E. Kamasak.
\newblock {A directional total variation}.
\newblock {\em Eur. Signal Process. Conf.}, 19(12):265--269, 2012.

\bibitem{Beck2009a}
A.~Beck and M.~Teboulle.
\newblock {A fast iterative shrinkage-thresholding algorithm for linear inverse
  problems}.
\newblock {\em SIAM J. Imaging Sci.}, 2(1):183--202, jan 2009.

\bibitem{Berkels2006}
B.~Berkels, M.~Burger, M.~Droske, O.~Nemitz, and M.~Rumpf.
\newblock {Cartoon extraction based on anisotropic image classification}.
\newblock {\em Proc. Vision, modeling and visualization}, pages 293--300, 2006.

\bibitem{Boyd2010}
S.~Boyd, N.~Parikh, E.~Chu, B.~Peleato, and J.~Eckstein.
\newblock {Distributed optimization and statistical learning via the
  alternating direction method of multipliers}.
\newblock {\em Found. Trends{\textregistered} in Mach. Learn.}, 3(1):1--122,
  2010.

\bibitem{Bredies2014}
K.~Bredies.
\newblock {\em Recovering piecewise smooth multichannel images by minimization
  of convex functionals with total generalized variation penalty}, pages
  44--77.
\newblock Springer Berlin Heidelberg, Berlin, Heidelberg, 2014.

\bibitem{Bredies2010}
K.~Bredies, K.~Kunisch, and T.~Pock.
\newblock {Total generalized variation}.
\newblock {\em SIAM J. Imaging Sci.}, 3(3):492--526, 2010.

\bibitem{Bredies2011}
K.~Bredies and T.~Valkonen.
\newblock {Inverse problems with second-order total generalized variation
  constraints}.
\newblock {\em Int. Conf. Sampl. Theory Appl.}, pages 1--4, 2011.

\bibitem{Chambolle2011a}
A.~Chambolle and T.~Pock.
\newblock {A first-order primal-dual algorithm for convex problems with
  applications to imaging}.
\newblock {\em J. Math. Imaging Vis.}, 40(1):120--145, dec 2011.

\bibitem{Chan2000}
T.~Chan, A.~Marquina, and P.~Mulet.
\newblock {High-order total variation-based image restoration}.
\newblock {\em SIAM J. Sci. Comput.}, 22(2):503--516, 2000.

\bibitem{Chan1998}
T.~Chan and C.~K. Wong.
\newblock {Total variation blind deconvolution.}
\newblock {\em IEEE Trans. Image Process.}, 7(3):370--5, 1998.

\bibitem{Delaney1998}
A.~H. Delaney and Y.~Bresler.
\newblock {Globally convergent edge-preserving regularized reconstruction: An
  application to limited-angle tomography}.
\newblock {\em IEEE Trans. Image Process.}, 7(2):204--221, 1998.

\bibitem{dong2009efficient}
Y.~Dong, M.~Hinterm{\"u}ller, and M.~Neri.
\newblock An efficient primal-dual method for ${L}^1${TV} image restoration.
\newblock {\em SIAM J. Appl. Math.}, 2(4):1168--1189, 2009.

\bibitem{dong2013convex}
Y.~Dong and T.~Zeng.
\newblock A convex variational model for restoring blurred images with
  multiplicative noise.
\newblock {\em SIAM J. Appl. Math.}, 6(3):1598--1625, 2013.

\bibitem{Easley2009}
G.~R. Easley, D.~Labate, and F.~Colonna.
\newblock {Shearlet-based total variation diffusion for denoising}.
\newblock {\em IEEE Trans. Image Process.}, 18(2):260--268, 2009.

\bibitem{Esedoglu2004}
S.~Esedoglu and S.~J. Osher.
\newblock {Decomposition of images by the anisotropic Rudin-Osher-Fatemi
  model}.
\newblock {\em Commun. Pure Appl. Math.}, 57(12):1609--1626, 2004.

\bibitem{Estellers2015}
V.~Estellers, S.~Soatto, and X.~Bresson.
\newblock {Adaptive regularization with the structure tensor}.
\newblock {\em IEEE Trans. Image Process.}, 24(6):1777--1790, 2015.

\bibitem{Fei2013}
X.~Fei, Z.~Wei, and L.~Xiao.
\newblock {Iterative directional total variation refinement for compressive
  sensing image reconstruction}.
\newblock {\em Signal Process. Lett. IEEE}, 20(11):1070--1073, 2013.

\bibitem{Fernandez-Granda2013}
C.~Fernandez-Granda and E.~J. Candes.
\newblock {Super-resolution via transform-invariant group-sparse
  regularization}.
\newblock {\em 2013 IEEE Int. Conf. Comput. Vis.}, pages 3336--3343, 2013.

\bibitem{Ferstl2013}
D.~Ferstl, C.~Reinbacher, R.~Ranftl, M.~Ruether, and H.~Bischof.
\newblock {Image guided depth upsampling using anisotropic total generalized
  variation}.
\newblock {\em Proc. IEEE Int. Conf. Comput. Vis.}, pages 993--1000, 2013.

\bibitem{Granlund1995}
G.~H. Granlund and H.~Knutsson.
\newblock {\em {Signal Processing for Computer Vision}}.
\newblock Springer-Science+Business Media Dordrecht, 1995.

\bibitem{Hafner2015}
D.~Hafner, C.~Schroers, and J.~Weickert.
\newblock {Introducing maximal anisotropy into second order coupling models}.
\newblock {\em Ger. Conf. Pattern Recognit.}, 9358:79--90, 2015.

\bibitem{Holler2014}
M.~Holler and K.~Kunisch.
\newblock {On infimal convolution of total variation type functionals and
  applications}.
\newblock {\em SIAM J. Imaging Sci. J. Imaging Sci.}, 7(4):2258--2300, 2014.

\bibitem{Jespersen2016}
K.~M. Jespersen, J.~Zangenberg, T.~Lowe, P.~J. Withers, and L.~P. Mikkelsen.
\newblock {Fatigue damage assessment of uni-directional non-crimp fabric
  reinforced polyester composite using X-ray computed tomography}.
\newblock {\em Compos. Sci. Technol.}, 136:94--103, 2016.

\bibitem{Jespersen2016b}
K.~M. Jespersen, J.~Zangenberg, T.~Lowe, P.~J. Withers, and L.~P. Mikkelsen.
\newblock {X-ray CT Data: Fatigue Damage in Glass Fibre/polyester Composite
  Used for Wind Turbine Blades [Data-set]}, 2016.

\bibitem{Jonsson}
E.~Jonsson, T.~Chan, and S.-C. Huang.
\newblock {Total variation regularization in positron emission tomography}.
\newblock Technical report, Dept. Mathematics, Univ. California., Los Angeles.,
  1998.

\bibitem{Lefkimmiatis2015}
S.~Lefkimmiatis, A.~Roussos, P.~Maragos, and M.~Unser.
\newblock {Structure tensor total variation}.
\newblock {\em SIAM J. Imaging Sci.}, 8(2):1090--1122, 2015.

\bibitem{Nesterov1983}
Y.~Nesterov.
\newblock {A method of solving a convex programming problem with convergence
  rate O (1/k2)}.
\newblock {\em Sov. Math. Dokl.}, 27(2):372--376, 1983.

\bibitem{Nikolova2000}
M.~Nikolova.
\newblock {Local strong homogeneity of a regularized estimator}.
\newblock {\em SIAM J. Appl. Math.}, 61(2):633--658, 2000.

\bibitem{Ranftl2012}
R.~Ranftl, S.~Gehrig, T.~Pock, and H.~Bischof.
\newblock {Pushing the limits of stereo using variational stereo estimation}.
\newblock {\em IEEE Intell. Veh. Symp. Proc.}, (1):401--407, 2012.

\bibitem{Ring2000}
W.~Ring.
\newblock Structural properties of solutions to total variation regularization
  problems.
\newblock {\em ESAIM: M2AN}, 34(4):799--810, 2000.

\bibitem{Rudin1992}
L.~I. Rudin, S.~Osher, and E.~Fatemi.
\newblock {Nonlinear total variation based noise removal algorithms}.
\newblock {\em Phys. D Nonlinear Phenom.}, 60(1-4):259--268, 1992.

\bibitem{Sandoghchi2014}
S.~R. Sandoghchi, G.~T. Jasion, N.~V. Wheeler, S.~Jain, Z.~Lian, J.~P. Wooler,
  R.~P. Boardman, N.~K. Baddela, Y.~Chen, J.~R. Hayes, E.~N. Fokoua,
  T.~Bradley, D.~R. Gray, S.~M. Mousavi, M.~N. Petrovich, F.~Poletti, and D.~J.
  Richardson.
\newblock {X-ray tomography for structural analysis of microstructured and
  multimaterial optical fibers and preforms}.
\newblock {\em Opt. Express}, 22(21):26181, 2014.

\bibitem{Scherzer1998}
O.~Scherzer.
\newblock {Denoising with higher order derivatives of bounded variation and an
  application to parameter estimation}.
\newblock {\em Computing}, 60(1):1--27, 1998.

\bibitem{sciacchitano2015variational}
F.~Sciacchitano, Y.~Dong, and T.~Zeng.
\newblock Variational approach for restoring blurred images with {C}auchy
  noise.
\newblock {\em SIAM J. Appl. Math.}, 8(3):1894--1922, 2015.

\bibitem{Setzer2008b}
S.~Setzer and G.~Steidl.
\newblock {Variational methods with higher order derivatives in image
  processing}.
\newblock {\em Approx. Theory XII San Antonio 2007}, pages 360--385, 2008.

\bibitem{Setzer2008}
S.~Setzer, G.~Steidl, and T.~Teuber.
\newblock {Restoration of images with rotated shapes}.
\newblock {\em Numer. Algorithms}, 48(1-3):49--66, 2008.

\bibitem{Setzer2011}
S.~Setzer, G.~Steidl, and T.~Teuber.
\newblock {Infimal convolution regularizations with discrete ll-type
  functionals}.
\newblock {\em Commun. Math. Sci.}, 9(3):797--827, 2011.

\bibitem{Temam1985}
R.~Temam.
\newblock {\em {Mathematical Problems in Plasticity}}.
\newblock Gaulthier-Villars, 1985.

\bibitem{Turgay2009}
E.~Turgay and G.~B. Akar.
\newblock {Directionally adaptive super-resolution}.
\newblock {\em 2009 16th IEEE Int. Conf. Image Process.}, (1):1201--1204, 2009.

\bibitem{Vogel1996}
C.~R. Vogel and M.~E. Oman.
\newblock {Iterative methods for total variation denoising}.
\newblock {\em SIAM J. Sci. Comput.}, 17(1):227--238, 1996.

\bibitem{Weickert1998a}
J.~Weickert.
\newblock {Anisotropic diffusion in image processing}.
\newblock {\em Image Rochester NY}, 256(3):170, 1998.

\bibitem{Weickert1999a}
J.~Weickert.
\newblock {Coherence-enhancing diffusion filtering}.
\newblock {\em Int. J. Comput. Vis.}, 31(2):111--127, 1999.

\bibitem{Weickert1998}
J.~Weickert, B.~M. T.~H. Romeny, and M.~A. Viergever.
\newblock Efficient and reliable schemes for nonlinear diffusion filtering.
\newblock {\em IEEE Trans. Image Process.}, 7:398--410, 1998.

\bibitem{Zhang2013}
H.~Zhang and Y.~Wang.
\newblock {Edge adaptive directional total variation}.
\newblock {\em J. Eng.}, (October):1--2, 2013.

\end{thebibliography}

\end{document}